\documentclass[10pt]{amsart}
\usepackage{amssymb}
\usepackage{amsmath}
\usepackage{amsthm}
\usepackage{epsfig}
\usepackage{mathrsfs}
\usepackage{graphicx}
\usepackage{ textcomp }
\usepackage{xcolor}
\usepackage{hyperref,comment}
\usepackage{bbold}
\usepackage{tikz}
\usepackage{amssymb}

\numberwithin{equation}{section}

\newtheorem{theorem}{Theorem}[section]
\newtheorem{theorem*}{Theorem}
\newtheorem{remark}[theorem]{Remark}
\newtheorem{lemma}[theorem]{Lemma}

\newtheorem{proposition}[theorem]{Proposition}

\newtheorem{definition}[theorem]{Definition}

\definecolor{darkgreen}{rgb}{0,0.4,0}

\def\supp{\mathrm{supp}\,}

\def\NN{{\mathbb N}}

\def\R{{\mathbb R}}

\def\eps{\varepsilon}

\def\e{\varepsilon}
\def\vphi{\varphi}

\def\na{\nabla}
\def\pa{\partial}

\begin{document}

\title{
Incompressible limit of a porous media equation with bistable and monostable reaction term
}
\author[I. Kim]{Inwon Kim}
\address{Department of Mathematics, UCLA,  Los Angeles, CA} 
\email{ikim@math.ucla.edu}

\author[A.Mellet]{Antoine Mellet}
\address{Department of Mathematics, University of  Maryland, College Park, MD}
\email{mellet@umd.edu}

\thanks{I. Kim was partially supported by NSF Grant DMS-1900804. \\
A. Mellet was partially supported by NSF Grant DMS-2009236.
}

\maketitle
\begin{abstract}
We study the incompressible limit of the porous medium equation with a reaction term that is non-monotone with respect to the pressure variable. More specifically we consider reaction terms that are either bistable or monostable. We show that this type of reaction term generates many interesting differences in the qualitative behavior of solutions, in contrast to the problem with monotone reaction terms that have been extensively studied in recent literature.  
After characterizing the limit problem, we embark on a comprehensive study of the problem in one space dimension, to illustrate the delicate nature of the problem, including the generic nature of non-uniqueness and instability. 
For compactly supported initial data, we show that the density can either perish or thrive, even if it starts from the same initial data, depending on its initial pressure configuration. 
When the initial pressure is a characteristic function, we establish the existence of the sharp threshold separating the two behaviors. Lastly we present a detailed analysis of the behavior of traveling waves in this incompressible limit. We study the existence of traveling waves for the limiting model and prove convergence results  in the incompressible limit
(depending on the reaction term).

 \end{abstract}

\section{Introduction} 

\subsection{A degenerate reaction diffusion equation and its incompressible limit}
We consider the following porous media equation with a reaction term:
\begin{equation}\label{rho_m}
\begin{cases}
\pa_t \rho - \mathrm{div}\, (\rho\na p)=  \rho f(p), \quad p =P_m(\rho) := \frac{m}{m-1}\rho^{m-1} & \mbox{ in } \R^n\times\R_+\\
\rho(x,0) =\rho_{in,m}(x) & \mbox{ in } \R^n
\end{cases}
\end{equation} 
where $f$ is a Lipschitz continuous function.
This is a very classical equation which arises in many applications, and it has been studied extensively in recent years as a mechanical model for tumor growth, see for instance \cite{PQV}, \cite{PV} and references therein. In this context several papers have investigated the asymptotic behavior of the solution of \eqref{rho_m} in the so-called incompressible regime, corresponding to the limit $m \to \infty$
 \cite{AKY,DP,DS,GKM,KP,PQV,PQTV,DHV} 
In the context of tumor growth, it is typically assumed that $p\mapsto f(p)$ is monotone  decreasing and satisfies
\begin{equation}\label{eq:fpm}
 f(p) \leq 0\quad  \mbox{  for }  p \geq p_M, 
\end{equation}
Under such a monotonicity assumption, the nonlinearity $g(\rho ) := \rho  f(p) = \rho f(\frac{m}{m-1}\rho^{m-1} )$ is a Fisher-KPP nonlinearity \cite{AW,KPP}. Indeed, we trivially have $g(\rho)>0$ in $(0, \rho_M)$, $g(\rho)<0$ in $(-\infty,0)\cup (\rho_M,+\infty)$ and the condition $g(\rho) \leq g'(0)\rho$ is equivalent to
$f(p)\leq f(0)$ for all $p\geq 0$.
\medskip

The goal of this paper is to understand the dynamic of this problem, still in the regime $m\gg1$, but when we {\bf remove} the monotonicity assumption on $f$ (but still assume \eqref{eq:fpm} to have uniform bounds on the pressure).
A simple example that we have in mind is the following quadratic function:
\begin{equation}\label{eq:bistable}
 f(p) = (1-p)(p-\alpha), \quad \alpha\in [0,1), \quad \int_0^1 f(s)\, ds>0,
 \end{equation}
 although our analysis will not be restricted to this case.
The reaction-diffusion equation \eqref{rho_m}-\eqref{eq:bistable} can be used to model the propagation of a biological population whose tendency to disperse depends on the population density when the nonlinearity \eqref{eq:bistable} takes into account a  weak ($\alpha=0$) or strong ($\alpha>0$) Allee effect (which assumes a reduced, or negative, growth rate at low density \cite{DK}).
When $ \alpha\in(0,1)$, 
 the function $g(\rho)= \rho f(p)=\rho f\left( \frac{m}{m-1}\rho^{m-1}\right) $ is a {\bf bistable}  reaction term while when   $\alpha=0$, $g$ is a {\bf monostable}  nonlinearity. 

In the context of tumor growth, such a function $f$ can be used to model the competition between drugs invading the tumor and the contact inhibition of tumor cells. 
Indeed, in classical tumor models, drugs are injected far from the tumor and diffuse toward the tumor (see \cite{BL}). They therefore act mostly on the cells in the regions with lower pressure, namely near the periphery of the tumor, while deep inside the tumor the natural growth due to cell division dominates the dynamic. Naturally one could ask whether this problem yields a stable evolution of the tumor zone. We will prove that it is not the case and that non-uniqueness may occur in the incompressible limit. Furthermore we will see that the tumor region may disappear or thrive with the same initial configuration, depending on its initial pressure. This is in contrast to the case with monotone $f$, where the evolution of the tumor is uniquely determined by its initial configuration.
\medskip

Note that we could easily extend our result to more general settings, in particular to inhomogeneous $f$ that depends also on time and space.  In the context of tumor growth, this could include models that take into account spatial in-homogeneities or dependence of $f$ on the concentration of nutrient. We will keep the framework simple to focus on the non-monotonicity of $f$ with respect  to the pressure variable and its consequences.

\medskip

\subsection{Known results when $p\mapsto f(p)$ is  decreasing}

Before stating our main results, we briefly recall here the most important aspects of the analysis  when $p\mapsto f(p)$ is monotone decreasing (and satisfies  \eqref{eq:fpm}): 
For appropriate initial data $\rho_{in,m}$, satisfying in particular $p_{in,m}\leq p_M$ and  $\rho_{in,m}\to \rho_{in}$,
the following basic properties have been proved 
\begin{enumerate}
\item The functions $\rho_m(x,t)$ and $p_m(x,t)$ are bounded in $L^\infty(\R^n\times\R_+)$ and in $BV(\R^n\times(0,T))$ and converge (strongly in $L^1_{loc} $ and a.e.)  to $\rho_{\infty}(x,t)$ and $p_\infty(x,t)$.
\item The pair $(\rho_{\infty},p_\infty)$ is solution (in the sense of distribution) of 
\begin{equation}\label{eq:asf}
\begin{cases}
\pa_t \rho_\infty = \Delta p_\infty + \rho_\infty f(p_\infty) & \mbox{ in } \R^n\times\R_+ \\
\rho_\infty (x,0)=\rho_{in}(x)& \mbox{ in } \R^n
\end{cases}
\end{equation}
and satisfies the Hele-Shaw graph condition
\begin{equation}\label{eq:HSG} 
p_\infty\in P_\infty(\rho_\infty)
:= \left\{ 
\begin{array}{ll} 
0 & \mbox{ if } \rho_\infty<1 \\ 
{[}0,\infty) & \mbox{ if } \rho_\infty=1  \\
\infty & \mbox{ if } \rho_\infty>1.
\end{array} 
\right.
\end{equation}
Furthermore, the problem \eqref{eq:asf}-\eqref{eq:HSG} has a unique weak solution.
\item For all $t>0$, $p_\infty(\cdot,t^+)$ (defined as the trace of the $BV$ function $p_\infty$) is the unique solution $p^*$ of the obstacle problem
\begin{equation}\label{eq:obstacle}
\begin{cases}
 p^*\geq 0, \quad \Delta p^* + f(p^*) \leq 0, \quad \mbox{ in } \{ \rho_{\infty}(\cdot,t)=1\}\\
 \Delta p^* + f(p^*) = 0 \mbox{ in } \{ p^*>0\}.
\end{cases}
 \end{equation}
\end{enumerate}
The first two points are proved in \cite{PQV}, while the characterization of $p_\infty(\cdot,t^+)$  as the solution of an elliptic obstacle problem \eqref{eq:obstacle} for all $t>0$ is proved in \cite{GKM}.
We can interpret \eqref{eq:asf}-\eqref{eq:HSG} as a Hele-Shaw type free boundary problem in which the 
saturated set $\Omega(t)=\{ \rho_{\infty}(\cdot,t)=1\}$ expands with normal velocity proportional to the pressure gradient  (see \cite{GKM} for a rigorous statement). 
We note that these results hold if $f$ depends on $(x,t)$ as well, provided the function $p\mapsto f(x,t,p)$ satisfies \eqref{eq:fpm} for all $(x,t)$. But the monotonicity of $p\mapsto f(p)$ plays a crucial role in the proof of the results above. 
Without it, we will see that neither the obstacle problem \eqref{eq:obstacle} nor the asymptotic equation \eqref{eq:asf}-\eqref{eq:HSG} can be expected to have a unique solution. 
Even the derivation of the $BV$ bounds, crucial to derive strong convergence of the density variable, appears to require this assumption.

\subsection{Main results: General case}
We will first state the result that one can show for general function $f$  satisfying the condition \eqref{eq:fpm} (but no monotonicity).
We start with the following proposition:
\begin{proposition}\label{prop:0}
Assume that $f$ is a Lipschitz function satisfying \eqref{eq:fpm} and that the initial condition $\rho_{in,m}$ is a non-negative function such that
\begin{equation}\label{eq:init}
 \rho_{in,m}\in L^1(\R^n),\qquad  p_{in,m} = P_m(\rho_{in,m}) \leq p_M \mbox{ a.e. in } \R^n.
 \end{equation}
Then, the solution $(\rho_m,p_m)$ of \eqref{rho_m} satisfies the following:
\item[(i)] $\rho_m(x,t)$ and $p_m(x,t)$ are uniformly bounded in $L^{\infty}(\R^n\times \R_+)\cap L^\infty(0,T;L^1(\R^n))$ for all $T>0$. Up to a subsequence, $\rho_m$ and $p_m$ converge  weakly-$*$ in $L^\infty$ to $\rho_\infty$ and $p_\infty$.
\item[(ii)] $p_m$  is uniformly bounded in $L^2(0,T;H^1(\R^n))$ and $\rho_m$ is uniformly bounded in $C^{1/2}(0,T;H^{-1}(\R^n))$ for all $T>0$.
\item[(iii)] Any accumulation point $(\rho_\infty,p_\infty)$ of $\{(\rho_m,p_m)\}_{m\geq 1}$ satisfies the Hele-Shaw graph condition \eqref{eq:HSG} a.e. in $\R^n\times\R_+$. 
\end{proposition}
However, the lack of $BV$ estimates, or more precisely the lack of strong convergence and  almost everywhere convergence of $p_m$, is hindering further characterization of the limit. 
In particular, passing to the limit in the equation \eqref{rho_m} to show that $(\rho_\infty,p_\infty)$ solves \eqref{eq:asf} is very delicate since the nonlinear term $f(p_m)$ might not converge to $f(p)$.
Nevertheless, we can observe that $f(p_m)$ is bounded in $L^\infty(\R^n\times\R_+)\cap L ^2(0,T;H^1(\R^n))$. We can thus assume, up to another subsequence, that there exists $g\in L^\infty(\R^n\times\R_+)\cap L ^2(0,T;H^1(\R^n))$ such that 
\begin{equation}\label{eq:fg}
f(p_m) \rightharpoonup g\qquad  \mbox{ weak-* in }L^\infty(\R^n\times \R_+).
\end{equation}

With this weak limit we can still characterize the limit pressure: 
\begin{proposition}\label{prop:1}
Under the assumption of Proposition \ref{prop:0} and for a subsequence such that \eqref{eq:fg} holds,
the limit $(\rho_\infty,p_\infty)$  satisfies, in the sense of distribution,
\begin{equation}\label{eq:asg}
\begin{cases}
\pa_t \rho_\infty = \Delta p_\infty + \rho_\infty g, \quad  p_{\infty}(1-\rho_{\infty}) = 0  & \mbox{ in } \R^n\times\R_+ \\
\rho_\infty (x,0)=\rho_{in}(x)& \mbox{ in } \R^n.
\end{cases}
\end{equation}

Furthermore, for a.e. $t>0$ the function $q=p_\infty(\cdot,t) $ is in $H^1(\R^n)$ and satisfies
\begin{equation}\label{varg}
\int_{\R^n} \na  q \cdot \na \phi -  g \phi\, dx  = 0 
 \end{equation} 
for all $\phi \in H^1(\R^n)$ such that $\phi(x)(1-\rho_\infty(x,t))=0$ a.e. $x\in\R^n$.
\item Finally, if $f$ is concave, then $g\leq f(p_\infty)$ and if $f$ is convex then $g\geq f(p_\infty)$: in particular, if $f$ is linear, then $g=f(p_\infty)$.
\end{proposition}

\medskip

Since $p_m$ is bounded in $L ^2(0,T;H^1(\R^n))$, we see that the lack of strong convergence of $p_m$ is due to the possibility of oscillating behavior in time. While we do not yet know if such oscillations actually happen,  it is plausible given the non-uniqueness of the expected limiting pressure equation $-\Delta p = f(p)$ (see the discussion below).  
This issue can be avoided if we can guarantee that the solution of \eqref{rho_m} is monotone in time which is true for some particular initial data. Indeed, we have:
\begin{proposition}\label{prop:monotone}
Assume that the initial pressure $p_{in,m}$ satisfy
$$ \Delta p_{in,m} +f(p_{in,m})\chi_{\{p_{in,m}>0\}}\geq 0 \mbox{ in } \R^n$$
for $m\geq m_0$.
Then $p_m$ and  $\rho_m$ are non-decreasing in time, and, up to a subsequence, 
 $p_m$ converges to $p_\infty$ strongly in $L^1_{loc}$ and a.e. in $\R^n\times\R_+$. In particular, we can assume that
 $$g=f(p_\infty)$$
 in Proposition \ref{prop:1}. In other words, $(\rho_\infty,p_\infty)$ solves
\eqref{eq:asf}-\eqref{eq:HSG}.
\end{proposition}

\medskip

In the cases in which we are able to prove that $g=f(p_\infty)$ (e.g. for linear $f$ or for monotone-in-time solutions), it is natural to ask whether  \eqref{eq:asf}-\eqref{eq:HSG} actually identifies the limit $(\rho_\infty,p_\infty)$ uniquely. This would imply the convergence of the whole sequence $(\rho_m,p_m)$. Such a result, which was proved in \cite{PQV} with monotone decreasing $f$, does not hold in general. As we will prove below in the simple one dimensional framework, such uniqueness is false even for relatively simple  $f$, such as the quadratic function given by \eqref{eq:bistable}. 
This is not so surprising. Indeed, equation \eqref{eq:asg} typically identifies a unique pressure function $p_\infty(\cdot,t)$ as the solution of the elliptic problem \eqref{varg}.
But without the monotonicity of $f$ this problem is not well-posed. 
Indeed, given a set $\Omega$, it is well known \cite{A,AR,Lions} that the Dirichlet problem
\begin{equation}\label{eq:pressurelimit}
\begin{cases}
-\Delta q = f(q) & \mbox{ in } \Omega \\
q=0 & \mbox{ on } \pa\Omega
\end{cases}
\end{equation}
is well posed if $f$ is monotone decreasing, or if $\mathrm{Lip} f \leq \lambda_\Omega$ where $\lambda_\Omega$ denotes the first eigenvalue of the laplacian in $\Omega$. But in general, even if $f$ is linear, this boundary value problem might  have no solution or multiple solutions.

\subsection{Monostable and bistable reaction term in dimension $1$}

For any given $m>1$, equation \eqref{rho_m} is well posed, but as $m\to \infty$, the discussion above suggests it is asymptotically close to a potentially ill-posed problem. This suggests a relatively unstable dynamic when $m\gg1$. 
In order to illustrate the non-uniqueness of the limit problem and to better understand the dynamic of $(\rho_m,p_m)$ when $m\gg1$, we will now focus on the simpler one-dimensional framework.
Throughout this section, we assume that $f$ satisfies one of the following conditions:

\medskip

Either  
\begin{equation}\label{eq:f1}
f(p)>0 \mbox{ in } (\alpha,p_M), \quad  f(p) <0 \mbox{ in } [0,\alpha)\cup (p_M,\infty)  \hbox{ for some } \alpha\in (0, p_M),
 \end{equation}
or (satisfied with $f$ given by \eqref{eq:bistable} with $\alpha=0$)
\begin{equation}\label{eq:f1'}
f(0)=0, \quad f(p)>0 \mbox{ in } (0,p_M), \quad  f(p) <0 \mbox{ in }  (p_M,\infty).
 \end{equation}

Above conditions are satisfied with $f$ given by \eqref{eq:bistable}, \eqref{eq:f1} for $\alpha \in (0, p_M]$ or \eqref{eq:f1'} for $\alpha = 0$. It is easy to check that the reaction term $\rho\mapsto \rho f(p)=\rho f\left( \frac{m}{m-1}\rho^{m-1}\right) $ in \eqref{rho_m}   is a bistable non-linearity in case \eqref{eq:f1} and a monostable nonlinearity in case \eqref{eq:f1'}.
In the sequel we will thus refer to \eqref{eq:f1} and \eqref{eq:f1'} respectively as the bistable and the monostable case.

\medskip

For the standard bistable reaction-diffusion equation (with linear diffusion), the dynamic strongly depends on the sign of the integral of the reaction term. For the nonlinear diffusion equation \eqref{rho_m}, the relevant integral (see for example Theorem \ref{thm:TW0} below) is 
$$ \int_0^{\rho_M} \rho^m f\left(\frac{m}{m-1}\rho^{m-1}\right)\,d\rho  = \int_0^{p_M} \frac 1 m \left(\frac{m-1}{m} p \right)^{\frac{2}{m-1}} f(p) \, dp.$$
In the regime $m\gg1$, 
many results below will thus depend on the sign of the integral  $\int_0^{p_M} f(p)\, dp$.
We also introduce the constant $K$ defined by
\begin{equation}\label{eq:K}
K:= \sup_{p>0} \frac{f(p)}{p}<\infty
\end{equation}
(since $f(0)\leq0$, we have in particular $K \leq \mathrm{Lip} (f)$).
\medskip
\medskip

 \paragraph{The pressure equation.}
The first step is to study the semi-linear boundary value problem \eqref{eq:pressurelimit} which determines the limiting pressure. In one dimension, we have:
\begin{equation}\label{eq:ob1d}
\begin{cases}
 -u'' = f(u) ,\quad  u\geq 0  \quad\mbox{ in } (0,L)\\
u(0)=u(L)=0.
\end{cases}\end{equation}

\begin{figure}
\begin{center}
\begin{tabular}{ccc}
	\includegraphics[width=0.5\textwidth]{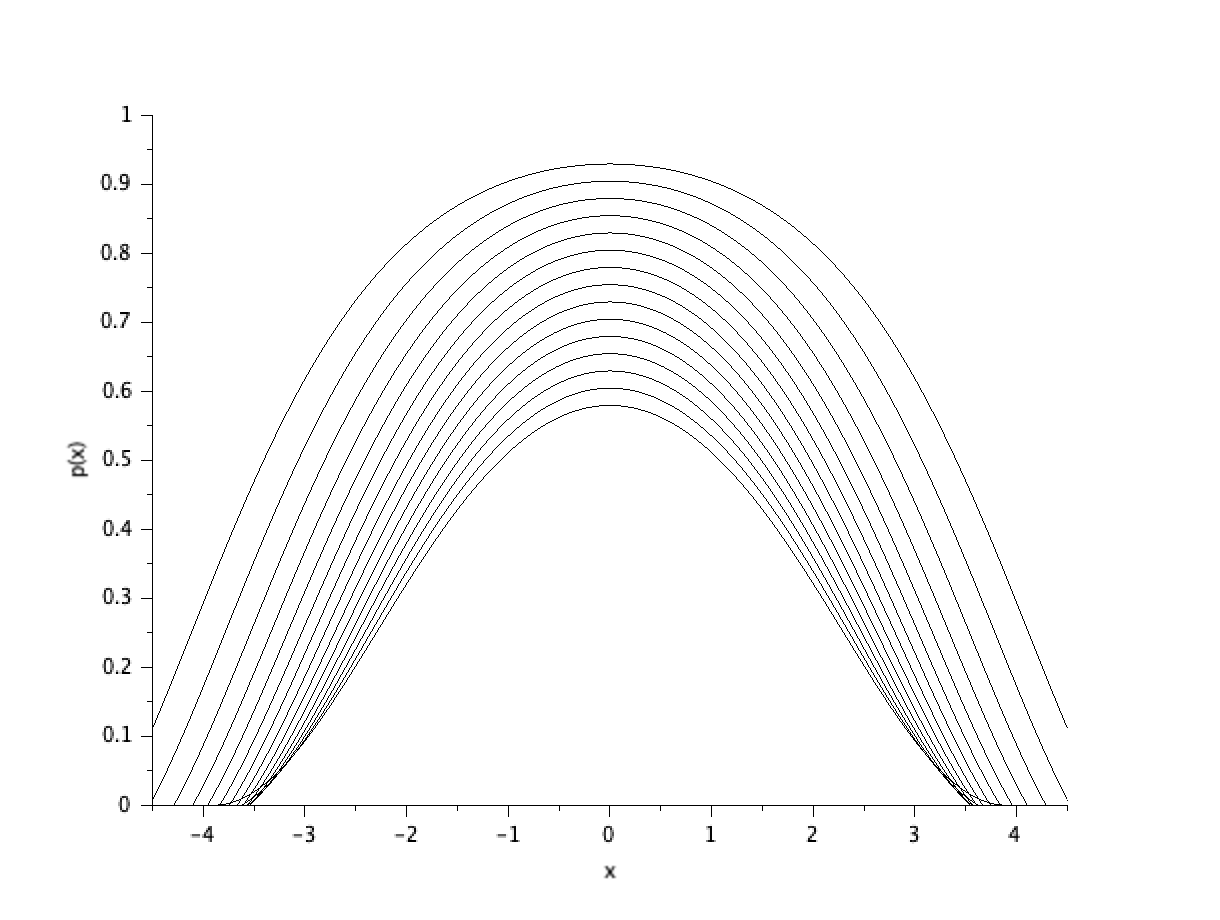} & 
	\includegraphics[width=0.2\textwidth]{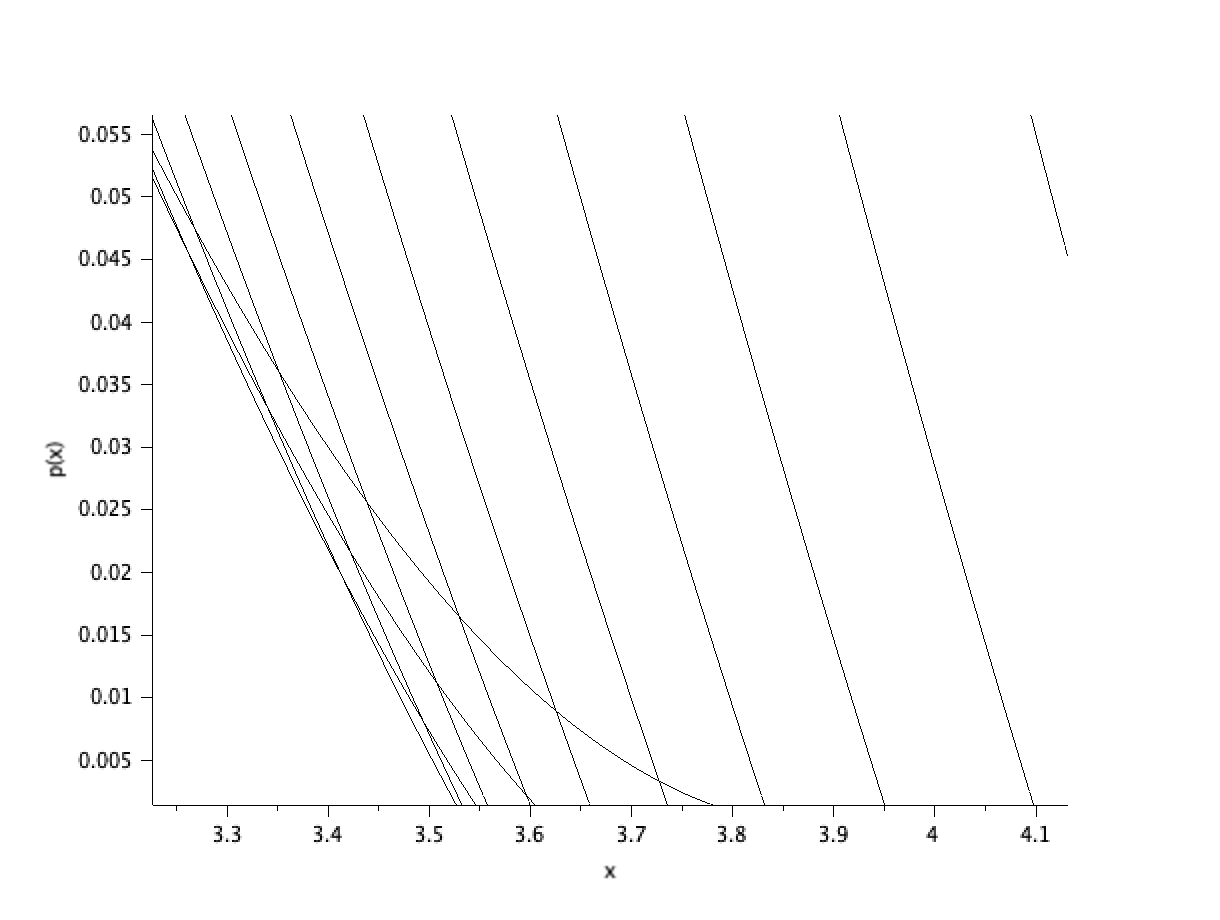}
\end{tabular}
	\caption{
	Solutions of equation \eqref{eq:ob1d} on $[-L,L]$ for various values of $L$ when $f(p)=(1-p)(p-1/4)$. The picture on the right, which is a close up near the boundary $x=L$, clearly shows that two solutions exists for some values of $L$.
		}\label{fig1}
\end{center}
\end{figure}

When $f$ satisfies \eqref{eq:f1'} ($\alpha=0$), then \eqref{eq:ob1d} as a solution $u=0$, which we referred to as the trivial solution below. Even in that case, it may also have non-trivial solution, as shown by the following result: 
\begin{proposition}\label{prop:L0}
Assume that $f$ is a Lipschitz function satisfying \eqref{eq:f1} or \eqref{eq:f1'}.
\item If  $ \int_0^{p_M} f(s)\, ds<0$ then \eqref{eq:ob1d} has no solution.
\item If  $ \int_0^{p_M} f(s)\, ds>0$, then there exists $L_0\geq \frac{\pi}{\sqrt K}$ (with $K$ given by \eqref{eq:K})
 such that
\begin{itemize}
\item[(i)] If $L<L_0$ then \eqref{eq:ob1d} has no  non-trivial solution.
\item[(ii)] If $L>L_0$ then \eqref{eq:ob1d} has at least one non-trivial  solution with support $(0,L)$.
\end{itemize}
\end{proposition}
Since  $u=0$ is only a solution of \eqref{eq:ob1d} in the monostable case \eqref{eq:f1'}
we see that  \eqref{eq:ob1d} has no solution in the bistable case  \eqref{eq:f1} when $L<L_0$.
We will also show that $L_0>\frac{\pi}{\sqrt K}$ except in the degenerate case where $f(p)=Kp$ for small $p>0$.
Numerical computations in the case $\alpha>0$ (see Figure \ref{fig1}) show that, at least for some $L>L_0$, \eqref{eq:ob1d}  might have two non-trivial solutions.  A rigorous verification of this non-uniqueness is given in Remark 4.2.

\medskip

This proposition proves that the characterization of $p_\infty$ given in Proposition \ref{prop:1} does not necessarily provide a unique characterization for the limit pressure.
In addition we will see below that, unlike the case studied in \cite{PQV,GKM} and others, the behavior of $(\rho_\infty,p_\infty)$ is not uniquely determined by the initial condition $\rho_{in}$. Note that even if we are given both  $\rho_{in}$ and $p_{in}$, it is not obvious that there is a unique solution for the limiting problem, but we do not address this issue here.

 \medskip
 
 \paragraph{Propagation vs extinction for compactly supported initial data.}
We establish this fact by studying the classical question of the long-term survival or extinction of bistable populations in this incompressible setting (in the framework of tumor growth, the question can be recast as: {\it who wins? the tumor or the drugs?}).

First we show that if the normalized initial pressure $p_{in,m}/p_M$ is a characteristic function, then a sharp threshold appears depending on the length of the support of $\rho_{in,m}$ when $f$ is in addition assumed to be concave (for instance when $f$ is given by \eqref{eq:bistable})
\begin{theorem}\label{thm:threshold}
Assume that $f$ is a Lipschitz function satisfying \eqref{eq:f1} or \eqref{eq:f1'}
 and that it is {\bf concave.}
Assume further that the initial condition is such that
\begin{equation}\label{eq:pinch} 
p_{in,m} = p_M\chi_{(0,L)} \qquad \mbox{ for some } L>0
\end{equation}
and let
$(\rho_\infty,p_\infty)$ be a limit point of $(\rho_m,p_m)$ with $L^\infty-*$ topology. Then
 the following holds:
\item[1.] If $ \int_0^{p_M} f(s)\, ds<0$, then $p_\infty(x,t)=0$ a.e. in $\R\times \R_+$ and 
$ \rho_\infty(t,x) = \rho_{in}(x) e^{f(0)t}$ for all $t>0$.
\item[2.] If $ \int_0^{p_M} f(s)\, ds>0$, then
\begin{itemize}
 \item[(i)] If $L<L_0$ then $p_\infty(x,t)=0$ a.e. in $\R\times \R_+$ and 
 $ \rho_\infty(t,x) = \rho_{in}(x) e^{f(0)t}$ for all $t>0$.
 \item[(ii)] If $L>L_0$, then    $\lim_{t\to\infty} \rho_\infty(t,x)=1$ for all $x\in\R$. 
\end{itemize}
\end{theorem}
We point out that when $ \int_0^{p_M} f(s)\, ds<0$ or $ \int_0^{p_M} f(s)\, ds>0$ and $L<L_0$, the limit $(\rho_\infty,p_\infty) = ( \rho_{in}(x) e^{f(0)t},0) $ is unique and therefore the whole sequence $(\rho_m,p_m)$ converges.
When $\alpha=0$, we then have $\rho_\infty(t) = \rho_{in}$ for all $t>0$, while when $\alpha>0$, we find
$\lim_{t\to\infty} \rho_\infty(t,x)=0$ for all $x\in\R$.
In the last case (ii), we can characterize the long time behavior of $\rho_\infty$ (for any accumulation point), but we do not prove that the limit is unique.

This sharp transition between extinction and propagation is reminiscent of a similar result for the reaction-diffusion equation with ignition non-linearity (see \cite{zlatos}).

\medskip

When \eqref{eq:pinch} holds,  the initial density satisfies
$$ \rho_{in,m} = \left(\frac{m}{m-1} p_M\right)^\frac{1}{m-1} \chi_{(0,L)} \to \chi_{(0,L)}.$$
However, we can have many different sequences of initial pressure $p_{in,m}$ (which do not satisfy \eqref{eq:pinch}) for which we also have $ \rho_{in,m}\to \chi_{(0,L)}$ and it is easy to see that the behavior of $\rho_\infty$ is not uniquely determined by $\rho_{in} = \lim \rho_{in,m}$, but also depends on $p_{in} = \lim p_{in,m}$. 
For example, when $ \int_0^{p_M} f(s)\, ds>0$ and  $L>L_0$, 
we can have \eqref{eq:pinch} and thus $\lim_{t\to\infty} \rho_\infty(t,x)=1$ as in Theorem \ref{thm:threshold}.
But we can also have $\rho_{in,m}\to \chi_{(0,L)}$ with $p_{in,m} (x) \leq \alpha-\delta$.
Since $\alpha-\delta$ is a supersolution for \eqref{rho_m}, we can easily show that any an accumulation point of $\{(\rho_m,p_m)\}_{m\geq 1}$ is such that $p_\infty(x,t)=0$  a.e. in $\R\times \R_+$ and $\lim_{t\to\infty} \rho_\infty(t,x)=0$ for all $x\in\R$.

This illustrates the unstable nature of the dynamic of \eqref{rho_m} when $m\gg1$.

\medskip

Finally, we note that Theorem \ref{thm:threshold} requires the assumption that $f$ is concave. When $f$ is not concave, it is not clear that we still have a sharp threshold, but we can prove the following:
\begin{proposition}\label{prop:long}
Assume that $f$ is a Lipschitz function satisfying \eqref{eq:f1} or \eqref{eq:f1'}
and that the initial condition satisfies \eqref{eq:init} and 
$$ p_{in,m} \leq p_M \chi_{(0,L)}\quad  \mbox{ for some } L<L^*:=\frac{\pi}{\sqrt K}$$
(with $K$ defined by \eqref{eq:K}). Then 
any accumulation point $(\rho_\infty,p_\infty)$ of $(\rho_m,p_m)$ (for the  $L^\infty-*$ topology) satisfies
$$p_\infty(x,t)=0\quad  \mbox{   and  } \quad \rho_\infty(x,t)=\rho_{in}(x) e^{f(0)t} \mbox{  a.e. in }\R\times \R_+.$$
\end{proposition}

\medskip

\paragraph{Traveling Wave Solutions.} 
Finally, we address a classical question with any reaction-diffusion equation: the existence and behavior of traveling wave solutions  joining the two equilibrium $\rho=0$ (and $p=0$) and $\rho_m^+ := \left(\frac{m-1}{m} p_M\right)^\frac{1}{m-1}$ (corresponding to $p=p_M$).
Traveling wave solutions of the reaction-diffusion equation \eqref{rho_m} are global in time solutions of the form  $\rho(x,t)=\bar \rho(x-ct)$ with $\bar \rho\in C(\R)$ satisfying
\begin{equation}\label{eq:eqTW} 
-c \bar \rho' = (\bar \rho \bar p')' + \bar \rho f(\bar p), \quad\mbox{ where } \bar p=\frac{m}{m-1}\bar \rho^{m-1} \qquad\mbox{ in } \R
\end{equation}
together with the limit conditions
\begin{equation}\label{eq:lim}
\begin{cases}
\lim_{y\to -\infty}\bar \rho (y) = \rho_m^+\\
\lim_{y\to + \infty} \bar \rho(y) = 0.
\end{cases}
\end{equation}

Existence and uniqueness of traveling waves for this degenerate diffusion-reaction equation
have been extensively studied, in particular, by Gilding and Kersner in \cite{Gilding_Kersner}. The results depend strongly on the nonlinearity $\rho\mapsto  \rho f\left(\frac{m}{m-1}\rho^{m-1}\right)$ and we summarize in the two theorems below the properties that are relevant to our analysis:
\begin{theorem}\label{thm:TW0}[Traveling wave of the $m$-equation - Bistable case]
 Assume that $f$ satisfies \eqref{eq:f1}.
For all $m\geq 1$, there exists a unique (up to translation)  solution $(\bar \rho_m,c^*_m)$ of \eqref{eq:eqTW}-\eqref{eq:lim}.  $\bar \rho_m$ is monotone decreasing  and the velocity $c^*_m$ has the same sign as
$$ \int_0^{\rho_m^+} \rho^m f\left(\frac{m}{m-1}\rho^{m-1}\right)\,d\rho .$$
Furthermore $\rho_m^+-\bar \rho_m(x)$ and its derivative up to order $2$ decay exponentially fast at $-\infty$ and 
\item[(i)] If $c^*_m<0$ then $\mathrm{Supp}\, \bar \rho_m = \R$ and $\bar \rho_m(x)$  and its derivative up to order $2$ decay exponentially fast at $+\infty$.
\item[(ii)] If $c^*_m\geq 0$ (and $m>1$), then (up to a translation) $\mathrm{Supp}\, \bar \rho_m = (-\infty,0)$ and when $c^*_m>0$, we have $\lim_{x\to 0^-} \bar p_m'(x) = -c^*_m$. 
\end{theorem} 
We refer to  \cite{Gilding_Kersner} for the proof of this result. In particular, 
the existence of a unique traveling wave in our framework follows  from Theorem 41 and Corollary 41.1. The 
 support properties of the traveling wave follow from  Theorem 45 and the limit $\lim_{x\to 0^-} \bar p_m'(x) = -c^*_m$ is proved in Theorem 46.
We point out that 
$$ 
 \int_0^{\rho_m^+} \rho^m f\left(\frac{m}{m-1}\rho^{m-1}\right)\,d\rho =
 \int_0^{p_M} \left( \frac{m-1}{m}p\right)^{\frac 2 {m-1}} f\left(p\right)\,\frac{1}{m}dp 
$$
which has the same sign, for $m$ large enough, as $\int_0^{p_M} f(p)\, dp$.

\begin{theorem}\label{thm:TW0'}[Traveling wave of the $m$-equation - Monostable case]
 Assume that $f$ satisfies \eqref{eq:f1'}.
For all $m\geq 1$ there exists $c^*_m$ such that \eqref{eq:eqTW}-\eqref{eq:lim} has a unique traveling wave with speed $c$  for all $c\geq c^*_m$, and no traveling waves exist when $c<c^*_m$.
For $c\geq c^*_m$, $ \bar \rho_m$ is monotone decreasing and $\rho_m^+-\bar \rho_m(x)$ and its derivative up to order $2$ decay exponentially fast at $-\infty$. Moreover
\item[(i)] If $c=c^*_m$, then 
$\mathrm{Supp}\,\bar  \rho_m = (-\infty,0)$ and $\lim_{x\to 0^-} \bar p_m'(x) = -c^*_m$.
\item[(ii)] If $c>c^*_m$, then $\mathrm{Supp}\, \bar \rho_m =\R$.
\end{theorem}

The existence of a unique traveling wave for $c\geq c^*$ follows from  
 \cite[Theorem 33]{Gilding_Kersner}. The support property of the traveling waves (points (i) and (ii)) follows from Theorem 38 in that same paper.

\medskip

Next, we consider the limit equation \eqref{eq:asf}-\eqref{eq:HSG}.
As before,  traveling waves are solutions of \eqref{eq:asf}-\eqref{eq:HSG}  of the form $\rho(x,t) =\bar \rho (x-ct)$. In particular, $\bar \rho$ must solve 
\begin{equation}\label{eq:limtw}
-c \bar \rho ' = \bar p'' + \bar \rho f(\bar p), \quad \bar p  \in P_\infty(\bar \rho),
\end{equation}
Furthermore, a traveling wave should connect equilibrium points, that is solutions of $ \bar \rho f(\bar p) = 0$, $\bar p  \in P_\infty(\bar \rho)$. 
We note that
$$ 
 \bar \rho f(\bar p) = \begin{cases} f(\bar p) & \mbox{ if } \bar \rho=1 \\ \bar \rho f(0) & \mbox{ if } \bar \rho<1
\end{cases}
$$
When $f$ satisfies \eqref{eq:f1} (bistable), the only solutions are $(\bar \rho ,\bar p)=(0,0)$ and $(\bar \rho ,\bar p) =(1,\alpha)$ or $(1,p_M)$. Since $(1,\alpha)$ is unstable, we will be looking for traveling waves connecting $(0,0)$ and $(1,p_M)$.
When $f$ satisfies \eqref{eq:f1'} (monostable), we have $f(0)=0$ so any $(\bar \rho ,\bar p)=(\ell,0)$ is an equilibrium point for all $\ell\in[0,1)$.
We will be looking for traveling waves connecting $(\ell,0)$ and $(1,p_M)$. Summarizing, we adopt the following definition:
\begin{definition}\label{def:TW}
We say that $(\bar \rho(x),\bar p (x))$ is a traveling wave of \eqref{eq:asf}-\eqref{eq:HSG} 
if 
$x\mapsto \bar p(x)$ and $x\mapsto \bar \rho(x)$ are non-increasing and there exists $c$ such that
\begin{equation}\label{eq:limtw}
-c \bar \rho ' = \bar p'' + \bar \rho f(\bar p), \quad \bar p  \in P_\infty(\bar \rho)\quad \hbox{ in the sense of distribution}
\end{equation}
with the boundary conditions
\begin{equation}\label{eq:TWbc}
 \begin{cases} \lim_{x\to -\infty} \bar \rho(x) = 1 \\ \lim_{x\to+\infty} \bar \rho(x)=\ell\end{cases}
 \mbox{ and }
\begin{cases} \lim_{x\to -\infty}\bar p (x) = p_M \\ \lim_{x\to+\infty} \bar p(x)=0\end{cases}  
\end{equation}
with $\ell=1$ in the bistable case and   $\ell \in [0,1)$ in the monostable case.
\end{definition}

\begin{remark}
It is likely unnecessary to assume that $\bar p$ and $\bar \rho$ are monotone to get uniqueness of the traveling waves. We choose to assume this given our context: since we are interested in the limits of the traveling waves given by Theorems \ref{thm:TW0} and \ref{thm:TW0'}, the monotonicity will be an immediate consequence of the monotonicity of $\bar \rho_m$.

Let us point out that the limiting problem has some interesting degeneracies  (see Proposition~\ref{prop:h}).
For example, if $f$  satisfies \eqref{eq:f1}  and $\int_0^{p_M} f(p) \, dp \geq 0$, then  there exists a solution of \eqref{eq:limtw} for all $c<0$  with periodic $\bar p(x)$ in $(-\infty,0)$ and monotone $\bar \rho$. This solution does not satisfy $\lim_{-\infty} \bar p = p_M$. Note also that when $f(0)=0$ (monostable case) and $c=0$, the equation  \eqref{eq:limtw} has infinitively many solutions found by taking $\bar p \equiv 0$ and any function $\bar \rho<1$.
\end{remark}

With this definition, we can prove the following proposition (to be compared with the classical result for the semilinear reaction-diffusion equation \cite{FM}):

\begin{proposition}\label{prop:limtw} 
We can classify the  travelling waves of  \eqref{eq:asf}-\eqref{eq:HSG}  as follows:
\begin{itemize}
\item[A.] Suppose $f$ is bistable, i.e. that it satisfies \eqref{eq:f1}. 
\begin{itemize}
\item[(i)] If $\int_0^{p_M} f(p) \, dp < 0$, then there exist no traveling waves of  \eqref{eq:asf}-\eqref{eq:HSG} (for any $c\in \R$).\\

\item[(ii)] If $\int_0^{p_M} f(p) \, dp = 0$, then there exist  a  traveling wave of  \eqref{eq:asf}-\eqref{eq:HSG} for all $c\leq 0$.
If $c<0$ then  $\mathrm{Supp}\, \bar \rho =\R$, while for $c=0$  we have $\mathrm{Supp}\, \bar \rho =(-\infty,0)$.
In both cases, $\ell=0$.


\item[(iii)] If $\int_0^{p_M} f(p) \, dp > 0$, then there exists a unique (up to translation) traveling wave of \eqref{eq:asf}-\eqref{eq:HSG}. It satisfies \eqref{eq:TWbc} with $\ell =0$ and
its speed is given by 
$$c^*:= \left(2\int_0^{p_M} f(p) \, dp \right)^{1/2}.$$
\end{itemize}
\item[B. ] Suppose $f$ is monostable, i.e. that it satisfies \eqref{eq:f1'}. Then for all $\ell \in[0,1)$
there exists
 a unique (up to translation)  traveling wave with non-zero speed satisfying \eqref{eq:TWbc}. Its speed is given by
$$c(\ell) = \frac{1}{1-\ell}  \left(2\int_0^{p_M} f(p) \, dp \right)^{1/2}\geq c^*.$$
\end{itemize}
\end{proposition}

We will also show that, up to translation, the traveling waves are given by $\bar \rho = \chi_{(-\infty,0)} + \ell \chi_{(0,+\infty)}$ and $\bar p$ as a solution $h$ of 
\begin{equation}\label{lim_tw}
 h''+f(h) = 0 \mbox{ in } (-\infty,0), \quad h(-\infty)=p_M, \quad h(0)=0.
 \end{equation}






\medskip

Lastly, we consider the traveling wave $(\bar \rho_m,c^*_m)$ of \eqref{rho_m} given either by Theorem \ref{thm:TW0} or by Theorem \ref{thm:TW0'} with the minimal speed. 
We then have the following result:
\begin{theorem}\label{thm:TW}

\item[(i)] If $f$  satisfies \eqref{eq:f1} and $\int_0^{p_M} f(p)\, dp<0$, then there exists $m_0$ and $\eta_0>0$ such that
\begin{equation}\label{eq:cm0}
 c^*_m \leq -\eta_0 m \qquad \mbox{ for all } m\geq m_0.
 \end{equation}
In particular $c^*_m\to -\infty$ and $p_m \to 0$ as $m\to\infty$.
\item[(ii)] If $f$  satisfies \eqref{eq:f1} and $\int_0^{p_M} f(p)\, dp> 0$, or if 
 $f$  satisfies \eqref{eq:f1'}, 
then 
$(\bar \rho_m,\bar p_m,c_m^*)$ converges to the unique traveling wave of \eqref{eq:asf}-\eqref{eq:HSG} with $\ell =0$. In particular
$$\lim_{m\to \infty} c^*_m =  \left(2\int_0^{p_M} f(p) \, dp \right)^{1/2}.$$ 
\end{theorem}
This theorem does not consider the case $\int_0^{p_M} f(p)\, dp= 0$. One issue in that case is that  the sign of $ \int_0^{\rho_m^+} \rho^m f\left(\frac{m}{m-1}\rho^{m-1}\right)\,d\rho $ might change as $m\to\infty$. Also missing in the statement is the behavior of the traveling waves with speed $c>c_m^*$ in the monostable case. We recall that for a given $c>c^*$, there exists a traveling wave $\bar \rho_m$ of \eqref{rho_m} (for $m$ large enough). The expectation is that, after translation, this traveling wave converges to  the traveling wave  of \eqref{eq:asf}-\eqref{eq:HSG}
with $\ell $ determined by the condition  $c =\frac{1}{1-\ell} c^*$ (that is $\ell = 1 -\frac{c^*}{c}$  given by Proposition~\ref{prop:limtw}-B.

\bigskip

{\bf Outline for the rest of the paper:}
In the next section we briefly recall and prove some well-known fundamental properties of the porous media equation \eqref{rho_m}. Section 3 is then devoted to the proof of our results in the general case (namely all the theorems and propositions presented in Section 1.3 above) while sections 4 and 5 focus on the particular one-dimensional case: In Section 4 we establish the asymptotic behavior of the solution for compactly supported initial data, while Section 5 is devoted to traveling waves solutions.

\medskip

\section{Preliminaries: The porous media equation \eqref{rho_m}}
We recall here some important properties of the porous media equation \eqref{rho_m} which we will use throughout the paper. We refer to \cite{Vazquez} for the proofs of these results (and much more).
First, the degenerate diffusion - reaction equation \eqref{rho_m} is well-posed. Indeed we have:
\begin{enumerate}
\item Let $\rho^k_{in}$ be a sequence of smooth approximations of $\rho_{in,m}$ in $L^1(\R^n)$, satisfying in particular $\rho^k_{in}\geq \frac 1 k$. Then  \eqref{rho_m} has a unique smooth solution $\rho^k$ with initial condition  $\rho^k_{in}$.
\item When $k\to\infty$, $\rho^k$ converge to a solution of  \eqref{rho_m} locally uniformly in $L^{\infty}(\R^n \times \R_+)$.
\item  The solution obtained this way is in fact the unique solution of \eqref{rho_m}.
\end{enumerate}
This last point follows from the following lemma which we prove below for the reader's sake:
\begin{lemma}\label{lem:compa}
Given $\rho_1(x,t)$ and $\rho_2(x,t)$ two solutions of \eqref{rho_m}, we have
$$
\int_{\R^n} (\rho_1 - \rho_2)_+ (x,t) dx \leq e^{C(m+1)t}\int_{\R^n} (\rho_1 - \rho_2)_+(x,0) dx \qquad \forall t>0,
$$
where $C$ depends on $\sup f$, $\sup f'_+$ and $\| \rho_i\|_{L^\infty}$ but  is independent of $m$. 
In particular, 
if $\rho_1(x,0) \leq \rho_2(x,0)$ then $\rho_1(x,t) \leq \rho_2(x,t)$ for all $t>0$.
\end{lemma}
This lemma implies in particular the uniqueness, comparison property and $L^1$ stability of \eqref{rho_m} for any $m>1$.
Of course, when $m\gg1 $, the stability estimate   blows up, which is  not surprising since we will prove that when $m\to\infty$ the limit equation exhibits non-uniqueness in some cases.
Nevertheless,  Lemma \ref{lem:compa} implies that if we  consider solutions with ordered initial data, they will be ordered in the limit $m\to\infty$.   
We also recall that stability results that are uniform in $m$ can be obtained when $p\mapsto f(p)$ is monotone decreasing (see \cite{PQV}). 

\begin{proof} 
The proof is classical (see \cite{Vazquez}): let $\beta_\delta(s)$ be a smooth approximation of $\mathrm{sign}^+(s)$  such that $\beta_\delta(s)=0$ if $s<0$, $0\leq \beta_\delta(s)\leq 1$ and $\beta_\delta'(s)\geq 0$ for all $s$. Denoting $w:= \rho_1^m - \rho_2^m$, we can write:
\begin{align*}
\int_{\R^n}  (\rho_1-\rho_2)_t \beta_\delta(w) dx  
&= \int_{\R^n} \Delta w \beta_\delta(w) dx + \int_{\R^n} (f(p_1)\rho_1 - f(p_2)\rho_2 ) \beta_\delta(w)  dx \\ 
 & =  -\int_{\R^n} |\nabla w|^2 \beta_\delta'(w) dx + \int_{\R^n} f(p_1) (\rho_1-\rho_2) \beta_\delta(w) dx + \int_{\R^n} (f(p_1)-f(p_2)) \rho_2 \beta_\delta(w) dx \\
 & \leq   (\sup f) \int_{\R^n}  (\rho_1-\rho_2)_+ dx +  \| \rho_2 \|_{L^{\infty}} (\sup f')  \int_{\R^n}  (p_1-p_2)_+ dx  \\ 
& \leq [ \sup f  + m\| \rho_2 \|_{L^{\infty}}\| \rho_1 \|^{m-2}_{L^{\infty}} (\sup f') ]\int_{\R^n} (\rho_1 - \rho_2)_+ dx.
\end{align*}
Passing to the limit $\delta \to0$ and using the fact that $\mathrm{sign}^+( \rho_1^m - \rho_2^m ) = \mathrm{sign}^+(\rho_1 - \rho_2)$, 
we deduce
$$\frac{d}{dt} \int_{\R^n}  (\rho_1-\rho_2)_+   dx  
 \leq C( m+1)\int_{\R^n} (\rho_1 - \rho_2)_+ dx,
$$
and the result follows.
\end{proof}

\section{General nonlinearity: Proofs of the main results}
\subsection{Proof of Propositions \ref{prop:0}}
We start with the following lemma, which  is a consequence of the comparison principle and Assumption \eqref{eq:fpm}:
\begin{lemma} \label{lem:3.1}
Under the assumptions of Proposition \ref{prop:0},
the functions
$\rho_m$ and $p_m$ are uniformly bounded in $L^{\infty}(\R_+\times\R^n)$ and $L^\infty((0,T);L^1(\R^n))$ with respect to $m$.  
\end{lemma}
\begin{proof}
We recall that the pressure $p_m$ solves
\begin{equation}\label{eq:pressure} 
\pa_t p = (m-1) p(\Delta p + f(p)) +|\na p|^2
\end{equation}
so the comparison principle and \eqref{eq:init} gives $0\leq p_m(x,t)\leq p_M$ and $0\leq \rho_m(x,t)\leq p_M^{\frac{1}{m-1}}$.
Integrating \eqref{rho_m} gives
$$\frac{d}{dt} \int_{\R^n} \rho_m(x,t)\, dx \leq \left ( \sup_{p\in[0,p_M]} f(p)\right)  \int_{\R^n} \rho_m(x,t)\, dx$$
which implies the bound for $\rho_m$ in $L^\infty((0,T);L^1(\R^n))$.
Since $p_m\leq \frac{m}{m-1} p_M^\frac{m-2}{m-1} \rho_m$, the corresponding bound for $p_m$ follows.

\end{proof}

\begin{lemma} \label{lem:comp}
Under the assumptions of Proposition \ref{prop:0}, the pressure $p_m$ satisfies (for all $T>0$):
\begin{equation}\label{pressure_l2}
\|\nabla p_m\|^2_{L^2(Q_T)} \leq \frac{m-1}{m-2} \left ( \sup_{p\in[0,p_M]} f(p)\right)  \|p_m\|_{L^1(Q_T)} + \frac{1}{m-2} \|p_{in,m}\|_{L^1(\R^n)}
\end{equation} 
and the density $\rho_m$ is uniformly bonded in $C^{1/2}(0,T; H^{-1}(\R^n))$.
\end{lemma}

\begin{proof}
Integrating the pressure equation \eqref{eq:pressure}  gives
$$
\frac{d}{dt} \int_{\R^n} p_m(x,t)\, dx = - (m-2) \int_{\R^n} |\na p_m|^2\, dx + (m-1) \int_{\R^n} p_m f(p_m)\, dx
$$
and the bound on $p_m$ follows by integrating with respect to $t\in [0,T]$.
Next, the density equation \eqref{rho_m} gives, for any $\varphi \in C^{\infty}_c(\R^n)$
$$
\frac{d}{dt}\int_{\R^n} \rho_m \varphi dx = \int_{\R^n} \rho_m \nabla p_m \cdot  \nabla \varphi dx + \int_{\R^n} f(p_m)\rho_m \varphi dx.
$$
It follows that for any $0\leq a <b \leq T$, we have
\begin{align*}
\left|\int_{\R^n} [\rho_m(\cdot,b)-\rho_m(\cdot,a)]\varphi dx\right| &\leq \|\rho_m\|_{L^{\infty}(Q_T)}\|\nabla p_m\|_{L^2(Q_T)} (b-a)^{1/2}\|\nabla \varphi\|_{L^2(\R^n)} \\ 
&\quad + (b-a)^{1/2}\| f(p_m)\rho_m^{1/2} \|_{L^{\infty}(Q_T)} \| \rho_m\|_{L^{1}(Q_T)}^{1/2}
 \|\varphi\|_{L^2(\R^n)}, \\ 
&\leq C |b-a|^{1/2} \|\varphi\|_{H^1(\R^n)},
\end{align*}
and the result follows.
\end{proof}

\medskip

Lemma \ref{lem:3.1} and \ref{lem:comp} immediately  imply the statements (i) and (ii) in Proposition \ref{prop:0}.
It only remains to prove (iii), namely the fact that any limit $(\rho_\infty,p_\infty)$ satisfies the {\bf Hele-Shaw graph condition} $$p_\infty \in P_\infty(\rho_\infty).$$

Classically this condition  is obtained by the relation $p_m(1-\rho_m)  = p_m - \left(\frac{m-1}{m}\right)^\frac{1}{m-1} p_m^{\frac{m}{m-1}} $, since the right hand side converges to zero as $m\to\infty$.
But  passing to the limit in the product $p_m \rho_m$ requires some extra care.
In \cite{PQV,GKM} this was easily done since  the $BV$ bounds provided some strong convergence.
We do not have such strong convergence here, but we note that
Lemma \ref{lem:comp} implies some time regularity for  $\rho_m$ and space regularity for $p_m$. This means that some space-time compensated compactness type argument could yield the convergence of $p_m \rho_m$ to $p_\infty \rho_\infty$ in $\mathcal D'$.
We present here a short proof which was first proposed in \cite{santambrogio}:

\begin{lemma}\label{lem:HS}
Under the conditions of Proposition \ref{prop:0}, 
any accumulation point $(\rho_\infty, p_\infty)$ of $(\rho_m,p_m)$ satisfies $0\leq \rho_\infty\leq 1$ and 
$$
\int_{\R^n}  p_\infty(x,t) (1-\rho_\infty(x,t)) dx = 0 \qquad \hbox{  for a.e. } t>0.
$$ 
\end{lemma} 

\begin{proof}
Throughout the proof, we fix a subsequence such that $(\rho_m,p_m)$ converge to $(\rho, p)$ weak-* in $L^\infty$.
Since $p_m\leq p_M$ 
we have
$$\liminf_{m\to\infty} \int_{\R^n} p_m(1-\rho_m) \, dx  \geq  \liminf_{m\to\infty}  \int_{\R^n} p_m\left(1-\left(\frac{m-1}{m}p_M \right)^{\frac{1}{m-1}}\right) \, dx \geq 0$$
and using the fact that $x^{\frac{m}{m-1}} \geq \frac{m}{m-1} x - \frac 1 {m-1}$, we can see that
$$\limsup_{m\to\infty}  \int_{\R^n} p_m(1-\rho_m) \, dx = \limsup_{m\to\infty}  \int_{\R^n}  \left[p_m - \left(\frac{m-1}{m}\right)^\frac{1}{m-1} p_m^{\frac{m}{m-1}} \right]\, dx \leq 0.
$$
Hence
$$
 \int_{\R^n} p_m(1-\rho_m) \, dx \to 0  \quad \mbox{ as $m\to \infty$}
$$
so the proof require us to show that $p_m\rho_m$ converges (at least in the sense of distribution) to $p \rho$. 
We now define the time-averaged functions
$$
p_m^{a,b}(x):= \frac{1}{b-a} \int _a^b p_m(x,t) dt\quad \hbox{ and } \quad p^{a,b}(x):= \frac{1}{b-a} \int_a^b p(x,t) dt
$$
and note that the bound  \eqref{pressure_l2} implies that $p_m^{a,b}$ converges to $p^{a,b}$ strongly in $L^2(\R^n)$. We then write
\begin{align*}
\frac{1}{a-b}\int_a^b \int_{\R^n} p_m(1-\rho_m) dx dt 
& = \frac{1}{b-a} \int_a^b \int_{\R^n} p_m(x,t) (1-\rho_m(x,a) )dx dt \\
&\quad + \frac{1}{b-a} \int_a^b \int_{\R^n} p_m(x,t)(\rho_m(x,a) - \rho_m(x,t) )dx dt.
\end{align*}
The first integral in the right side can also be written as  $\int p_m^{a,b}(x) (1-\rho_m(x,a) )dx $ and thus
 converges to $\int p^{a,b}(x)(1-\rho(x,a)) dx$ when $m\to\infty$.. As for the second integral, we use the uniform density bound in $C^{1/2}([0,T]; H^{-1}(\R^n))$:
\begin{align*}
\int_a^b \int_{\R^n} p_m(x,t)(\rho_m(x,a) - \rho_m(x,t)) dx \, dt 
& \leq  \int_a^b \|\nabla p_m(\cdot,t)\|_{L^2(\R^n)} \|\rho_m(\cdot,a)-\rho_m(\cdot,t)\|_{H^{-1}(\R^n)} dt \\ 
& \leq \left(\int_a^b \|\nabla p_m(\cdot,t)\|^2_{L^2(\R^n)} dt\right)^{1/2} \left(\int_a^b \|\rho_m(\cdot,a)-\rho_m(\cdot,t)\|_{H^{-1}(\R^n)} ^2\, dt\right)^{1/2} \\
& \leq C \left(\int_a^b \|\nabla p_m(\cdot,t)\|^2_{L^2(\R^n)} dt\right)^{1/2} \left(\int_a^b (t-a)\, dt\right)^{1/2} \\
& \leq C \left(\int_a^b \|\nabla p_m(\cdot,t)\|^2_{L^2(\R^n)} dt\right)^{1/2} (b-a)
\end{align*}
where $C$ is independent of $m$. 
Since $f_m(t):=\|\nabla p_m(\cdot,t)\|^2_{L^2(\R^n)}$ is uniformly bounded in $L^1([0,T])$, along a subsequence $f_m(t)$ weakly converges in $L^1$ to a measure $\mu$ such that $\mu([0,T]) <\infty$. Hence
\begin{align*}
\limsup_{m\to\infty} \frac{1}{b-a} \int_a^b \int_{\R^n} p_m(x,t)(\rho_m(x,a) - \rho_m(x,t)) dx dt 
& \leq C\lim_{m\to\infty} \left(\int_a^b \|\nabla p_m(\cdot,t)\|^2_{L^2(\R^n)} dt\right)^{1/2} \\
&\leq C (\mu[a,b])^{1/2}.
\end{align*}
We have thus proved:
\begin{align*}
0 & = \lim_{m\to\infty}\frac{1}{a-b}\int_a^b \int_{\R^n} p_m(1-\rho_m) dx dt\\
&  = \frac{1}{b-a} \int_{\R^n} p^{a,b}(x,t) (1-\rho(x,a) )dx dt + O(\mu[a,b])^{1/2}.
\end{align*}
Now we fix $a$ and send $b$ to $a$. Since $\int_0^T d\mu <\infty$, $\mu[a,b]$ tends to zero as $b\to a$ expect for  possibly countably many $a$. On the other hand we have 
$$
\int p^{a,b} (x) (1-\rho(x,a)) dx \to \int p(x,a)(1-\rho(x,a)) dx
$$
at Lebesgue points of $F(t):= \int p(x,t)(1-\rho(x,t)) dx$.
The result follows.
\end{proof}

\begin{remark}\label{rmk:1}
The proof  shows in particular that $\rho_m p_m$ converges to $\rho p$ in the sense of distribution.
We note that $|\na f(p_m)| = |f'(p_m) \na p_m| \leq \mathrm{Lip} f |\na p_m|$ and so $f(p_m)$ is bounded in $L^2(0,T;H^1(\R^n))$ just like $p_m$. 
So if we assume that 
$$ f(p_m) \rightharpoonup g\qquad  \mbox{ weak-* in }L^\infty( \R^n \times \R_+)$$
then the same argument shows that
$$ \rho_m f(p_m) \longrightarrow \rho g \quad \mbox{ in } \mathcal D'(\R^n \times \R_+).$$
\end{remark}

\subsection{Proof of Proposition \ref{prop:1}} 
Using Remark \ref{rmk:1} we can prove the following proposition, which is the first part of Proposition \ref{prop:1}:
\begin{proposition}\label{prop:weaklimit}
Under the assumptions of Proposition \ref{prop:0} and for a subsequence such that \eqref{eq:fg} holds,
the limit $(\rho_\infty,p_\infty)$ satisfies
\begin{equation}\label{eq:asg2}
\begin{cases}
\pa_t \rho_\infty = \Delta p_\infty + \rho_\infty g & \mbox{ in } \R^n\times\R_+ \\
\rho_\infty (x,0)=\rho_{in}(x)& \mbox{ in } \R^n
\end{cases}
\end{equation}
in the sense of distribution.
\end{proposition}
\begin{proof}
We can rewrite equation \eqref{rho_m} as
$$\pa_t \rho_m = \Delta \rho_m^m + \rho_m f(p_m)$$
where $\rho_m^m = \left( \frac{m-1}{m} p_m\right)^{\frac{m}{m-1}}$ converges to $p_\infty$ (since $p_m$ is bounded uniformly by $p_M$). Remark \ref{rmk:1} gives the convergence of $\rho_m f(p_m)$ to $\rho_\infty g$ and 
the result follows.
\end{proof}

The fact that $g\geq f(p_\infty)$ if $f$ is convex (and the opposite inequality when $f$ is concave) is classical:
We write $f(p) = \sup_{a\in \R} \{ a p -f^*(a)\}$ with $f^*$ the Legendre transform  of $f$.
This implies $ f(p_m(x,t)) \geq a p_m(x,t) -f^*(a)$ and so $g(x,t)\geq ap_\infty(x,t) -f^*(a)$ for all $a\in \R$. The result follows by taking the supremum over $a\in \R$.
In particular, when $f(p) = a p+b $ is linear we get that any accumulation point $(\rho_\infty,p_\infty)$ is a solution of 
$$
\pa_ t \rho =  \Delta p + \rho f(p), \quad \rho \leq 1,\,\,  \int p(\cdot,t)(1-\rho(\cdot,t)) dx =0  \hbox{ a.e. } t.
$$

\medskip

To complete the proof of Proposition~\ref{prop:1}, we must show that $p_\infty(\cdot,t)$  solves \eqref{varg} for a.e. $t>0$
With the lack of BV estimates on $\rho_{\infty}$ we cannot proceed as in \cite{GKM} to characterize pressure at all times as solutions of obstacle problems. But still it is possible to show that the pressure equation holds for a.e. time by proceeding as in \cite{MRS}:

\begin{proposition}\label{lem:concave}
Under the assumptions of Proposition \ref{prop:weaklimit},  for a e. $t>0$, $q=p_\infty(\cdot,t) $ solves 
$$
\int_{\R^n} \na  q \cdot \na \phi -  g \phi\, dx  = 0 
$$
for all $\phi \in H^1(\R^n)$ such that $\phi(x)(1-\rho_\infty(t,x))=0$ a.e. $x\in\R^n$.
In particular, $p_\infty$ satisfies
\begin{equation}\label{varg}
-\Delta p_\infty = g \qquad \hbox{ in }\quad  \mathrm{Int}(\{\rho_\infty(\cdot,t)=1\})  \quad \hbox{ for a.e. } t>0.
\end{equation} 
\end{proposition}

\begin{proof}
We write $(p,\rho)$ for $(p_{\infty}, \rho_{\infty})$ in this proof. For given $t_0>0$, let us define the functional space
$$H^1_{\rho(t_0)}:= \{\phi\in H^1(\R^n), \phi(1-\rho(\cdot,t_0))=0 \mbox{ a.e. in } \R^m\}.$$ 
\medskip

We choose $\phi \in H^1_{\rho(t_0)}$ such that $\phi\geq 0$.
Since $(\rho,p)$ solves $\rho_t -  \Delta p = g \rho$ in the weak sense with $p, g\in L^2([0,T]; H^1_0(\R^n))$ for any $T>0$, 
$$
\int_{t_0}^{t_0+h} \int \nabla \phi(x)\cdot  \nabla  p (x,t) - g(x,t)\rho(x,t)\, dx \, dt = \int [\rho(x,t_0) - \rho(x, t_0 + h)] \phi(x) dx \hbox{ for any } h>0.
$$
Since $\phi\in H^1_{\rho(t_0)}$, we have $\phi(x) \rho(x,t_0) = \phi(x)$ a.e.  and since  $\rho(x,t) \leq 1$ and $\phi(x)\geq 0$ a.e., we deduce
$$
\frac{1}{h}\int_{t_0}^{t_0+h} \int \nabla \phi(x)\cdot  \nabla  p (x,t) - g(x,t)\rho(x,t) \, dx \, dt  \leq 0 \hbox{ for any } h>0.
$$
Finally, we note that since $H^1(\R^n)$ is separable, Lebesgue's differentiation theorem implies that  a.e. $t_0>0$ is a differentiation point for the function
$ t\mapsto\int  \nabla \phi(x)\cdot  \nabla  p (x,t)$  for all $\phi\in H^1(\R^n)$.
We deduce that for a.e. $t_0>0$ (independent of $\phi$) we have
\begin{equation}\label{eq:phieq}
\int  \nabla \phi(x) \cdot \nabla p(x,t_0) dx \leq \int g(x,t_0)\rho(x,t_0) dx. 
\end{equation}
Hence we have $-\Delta p \leq g$ in  the interior of the set $\{\rho_\infty(\cdot,t)=1\})$  a.e. $t$.

\medskip

We prove the reverse inequality (still with $\phi\geq 0$) by arguing similarly by integrating from $t_0-h$ to $t_0$.
Once we have equality in \eqref{eq:phieq} for $\phi\geq 0$, it is easy to show that it holds for all $\phi$.

\end{proof}

\subsection{Proof of Proposition \ref{prop:monotone}}
For general $f$, it is still possible to recover $g=f(p_{\infty})$ for solutions with time monotonicity.
This is the result of Proposition \ref{prop:monotone} which we prove now:
\begin{proof}[Proof of Proposition \ref{prop:monotone}]
Suppose initial pressure $p_{in,m}$ satisfy
$$ \Delta p_{in,m} +f(p_{in,m})\chi_{\{p_{in,m}>0\}} \geq 0$$
for sufficiently large $m$. Then $h(x,t):=p_{in,m}(x)$ is a subsolution of the pressure equation \eqref{eq:pressure}, so the comparison principle yields $p_m(\cdot,\eps) \geq p_{in,m}(x) = p_m(\cdot,0)$ for any $\eps>0$. Using the comparison principle again, we deduce that for any $\e>0$ we have
$$
p_m(\cdot,t+\e) \geq p_m(\cdot,t) \hbox{ for } t\geq 0,
$$
or, in other words, $p_m$ (and thus $\rho_m)$ are non-decreasing functions. This monotonicity in time yields, together with the previous estimates obtained above, that $p_m$ is uniformly bounded in $W^{1,1}_{loc} (\R^n\times (0,T))$.
Indeed, since $\pa_t p_m \geq 0$ we can write
$$ 
\int_0^T \int_{\R^n} 
|\pa_t p_m | \, dx \, dt  =  \int_0^T \int_{\R^n} 
\pa_t p_m \, dx \, dt \leq \int_{\R^n} p_m(x,T)\, dx
 $$
and Lemma \ref{lem:3.1} implies that $\pa_t p_m$ (and similarly $\pa_t \rho_m$) is bounded in $L^1(\R^n\times (0,T))$.  Using Lemma~\ref{lem:comp}, we deduce that $p_m$ is uniformly bounded in $W^{1,1}_{loc} (\R^n\times (0,T))$.

Hence it follows that along subsequences we have strong convergence of $p_m$ in $L^1_{loc}(\R^n\times \R_+)$ and almost everywhere convergence.  In particular, we can assume that
 $g=f(p_\infty)$ in Proposition~\ref{prop:weaklimit}.  
 
 \end{proof}

\section{The bistable case in dimension $1$: Invasion vs. Extinction}
For the remainder of the paper, we assume that $n=1$ and that the nonlinearity $f$ satisfies \eqref{eq:f1} or \eqref{eq:f1'}.
We denote 
$$ F(p) = \int_0^p f(s)\, ds.$$

\subsection{Preliminaries: Proof of Proposition \ref{prop:L0}}
In this section, we investigate the properties of the elliptic equation  \eqref{eq:ob1d}, which we recall here for convenience:
\begin{equation}\label{eq:ob1d*}
\begin{cases}
 -u'' = f(u) ,\quad  u\geq 0  \quad\mbox{ in } (0,L)\\
u(0)=u(L)=0.
\end{cases}\end{equation}

\begin{proof}[Proof of Proposition \ref{prop:L0}]

We first make a couple of simple and classical remarks about \eqref{eq:ob1d*}:
Multiplying \eqref{eq:ob1d*} by $u'$ and integrating, we immediately see that 
any solution of \eqref{eq:ob1d*} satisfies:
\begin{equation}\label{eq:order1} 
\frac 1 2  |u'(x)|^2 + F(u(x)) = \frac{1}{2} |u'(0)|^2, \quad \forall x\in [0,L].
\end{equation}
If $u$ is a nontrivial solution of \eqref{eq:ob1d*}, there exists  $x_0$  such that  $u(x_0)=\sup_{x\in (0,L)} u(x)>0$. 
We have
$-u''(x_0)=f(u(x_0))\geq 0$ which gives in particular $u(x_0) \in [\alpha,p_M]$
and  $u'(x_0)=0$ which implies
$ \frac{1}{2} |u'(0)|^2 = F(u(x_0)) \leq F(p_M)$.
When $\alpha>0$, we cannot have $u(x_0)=\alpha$ since $F(\alpha)<0$ and if $u(x_0)= p_M$, then $u'(x_0)=u''(x_0)=0$ so \eqref{eq:ob1d*}  (and the fact that $f$ is Lipschitz) implies that $u\equiv p_M$, a contradiction.
We deduce 
$$u(x_0) \in (\alpha,p_M), \qquad u'(0)\in (0,\sqrt{2 F(p_M)}).$$

Finally, when $\alpha>0$, we observe that if $x_1\in(0,L)$ is such that $u(x_1)=0$, then \eqref{eq:ob1d*}  implies $-u''(x_1) = f(0)<0$ and so $x_1$ is a strict minimum of $u$. In particular, while we cannot claim that $u>0$ in $(0,L)$, we see that
the zeroes of $u$ are isolated in $(0,L)$ (when $\alpha=0$, the fact that $f$ is Lipschitz implies that either $u\equiv 0$ in $(0,L)$ or $u>0$ in $(0,L)$).

\medskip
When $F(p_M)<0$,  we have $F(p)<0$ for all $p\in(0,p_M]$ so \eqref{eq:order1} gives
$ \frac 1 2  |u'(x)|^2 \geq - F(u(x))>0$ for all $x$, which contradicts the existence of $x_0$ such that $u'(x_0)=0$. Hence no solution exists in that case.

\medskip

Next, we show that when $F(p_M)>0$, a non-trivial solution exists at least for some $L$: Given $ \gamma\in (0,\sqrt{2 F(p_M)})$,
let
$$G(s) : = \int_0^s  \frac{1}{ \sqrt{ \gamma^2 - 2 F(t)}} \, dt, \qquad s\in (0,s_0)
$$
where $s_0\in (\alpha,p_M)$ is such that $F(s_0) = \frac 1 2 \gamma^2$. We see that $G$ is a monotone increasing function such that
$$G(0)=0,  \quad G(s_0)=b_0\in(0,\infty), \quad G'(0) = \frac 1 \gamma,  \quad G'(s_0)=+\infty$$
(the fact that $b_0<\infty$ follows from the fact that $f(s_0)\neq 0$ and  
$\gamma^2 - 2 F(t) = 2 (F(s_0) -F(t) = 2 f(s_0)(s_0-t) + o(s_0-t)$).
and 
$$G'(s) > \frac {1}{\sqrt{2(F(p_M)-F(\alpha))}} \qquad \mbox{ for } s\in(0,s_0).$$
The function
$$u(x): = 
\begin{cases}
G^{-1}(x) & \mbox{ for }  x\in [0,b_0] \\ 
G^{-1}(2b_0-x) &  \mbox{ for }  x\in (b_0,2b_0]
\end{cases}
$$
then solves  \eqref{eq:ob1d*} with $L=2b_0$. 
\medskip
We now introduce the set
$$
S=\{L\, ;\,  \mbox{\eqref{eq:ob1d*} has at least one non-trivial solution in $(0,L)$}\}.
$$
The construction above shows that $S$ is non empty and Proposition \ref{prop:L0} will follow if we show that $S$ is an interval of the form $(L_0,\infty)$ or $[L_0,\infty)$ with $L_0>0$.

\medskip

Given $L_1\in S$ and $u_1$ solving  \eqref{eq:ob1d*}  in $(0,L_1)$, we construct a subsolution of \eqref{eq:ob1d*} in $(0,L)$ for $L>L_1$ by using the translation invariance of the equation.
First, we note that the extension of $u_1$ by $0$ (still denoted $u_1$) solves 
$$
-u_1''=f(u_1) \qquad \mbox{ in } \{u_1>0\}
$$
and so do its translations $u_1(x-y)$ for all $y\in\R$. We now define
\begin{equation}\label{superposition}
v_1 (x):= \sup\{ u_1(x-y)\, ;\, y \in [0,L-L_1]\} 
\end{equation}
which satisfies
$$ 
-v_1''\leq f(v_1) \qquad \mbox{ in } \{v_1>0\}
$$
and $ \{v_1>0\} = (0,L)$ because $u_1$ is non-negative in $(0,L_1)$ with isolated zeroes.

\medskip

Next, consider $v_2$ such that 
$$ -v_2'' = M, \mbox{ in } (0,L), v_2(0)=v_2(L)=0.$$
It is a supersolution for \eqref{eq:ob1d*} if $M\geq \sup f$ and satisfies $v_2 \geq v_1$ if $M$ is large enough. 
The existence of a non trivial solution of  \eqref{eq:ob1d*}  satisfying  $v_1\leq u\leq v_2$ now follows by Perron's principle.
\medskip

Finally, it only remains to show that $L_0=\inf \{L\, ;\, L \in S\}$ satisfies $L_0>0$.
For this, we assume that $u$ solves \eqref{eq:ob1d*} in $(0,L)$ for some $L>0$. Then
$$ \int_0^L |u'|^2\, dx = -\int_0^L u u''\,  dx= \int_0^L u f(u)\,  dx \leq K \int_0^L u^2  dx$$
with $ K= \sup_{p>0} \frac{f(p)}{p}<\infty$.
Poincar\'e inequality gives
$$ \int_0^L u^2  dx \leq  \left( \frac{L}{\pi}\right)^2
\int_0^L |u'|^2\, dx.$$
Together with the previous inequality, we get
$$
 \int_0^L |u'|^2\, dx  \leq \left( \frac{L}{\pi}\right)^2 K
\int_0^L |u'|^2\, dx.
$$  
If $u$ is not the trivial solution, we must have $  \left( \frac{L}{\pi}\right)^2K\geq 1$, that is
$L  \geq \frac{\pi}{\sqrt K}$.
Furthermore, if $  \left( \frac{L}{\pi}\right)^2K=1$ then the inequalities above must be equalities and so $f(u(x))= K(u(x))$ for all $x$ and so $f(p)=Kp$ for $p\in[0,\sup u]$.

\end{proof}


\begin{remark}\label{rmk:L_0} 
The non-existence parts of this result can be extended to  subsolutions, that is to functions satisfying 
$$
\begin{cases}
 -u'' \leq  f(u)  \quad\mbox{ in } (0,L)\\
u(0)=u(L)=0.
\end{cases}
$$
Indeed in that case \eqref{eq:order1} becomes
$\frac 1 2  |u'(x)|^2 \geq - F(u(x)) + \frac{1}{2} |u'(0)|^2$ which is always positive if $F(p_M)<0$, 
and this contradicts the existence of $x_0$ such that $u(x_0) = \sup_{(0,L)} u$.

\medskip

When $F(p_M)>0$ and $L<L_0$, if a subsolution $u$ exists in $(0,L)$ then we can use the supersolution $v_2$ as in the proof above (with $M$ very large so that $v_2>u$) and find a solution $w$ of \eqref{eq:ob1d*} in $(0,L)$ satisfying $u \leq w\leq v_2$, which contradiction the definition of $L_0$.
\end{remark}

\begin{remark}\label{rem:u0}
When $f$ is a bistable nonlinearity ($\alpha>0$), \eqref{eq:ob1d*} can have a solution with $u'(0)=u'(L)=0$. This solution can be found by solving 
$$ -u''=f(u),\quad  u(0)=0, \quad u'(0)=0.$$
Since $f$ is Lipschitz, this second order initial value problem has a unique solution, which satisfies
$\frac 1 2  |u'(x)|^2 + F(u(x)) =0$. It is not difficult to show that $u$ will be increasing until it reaches the value $\beta$ such that $F(\beta)=0$ and then decreasing until it goes back to zero and repeats that process periodically.
If we denote by $L_c$ the first time that $u$ reaches back to zero, we reduce that \eqref{eq:ob1d*} has solutions satisfying $u'(0)=u'(L)=0$ if and only if $L= k L_c$ for some $k\in \mathbb N$.  On the other hand, the construction given in the proof of Proposition~\ref{prop:L0} above provides  a solution $u$ with nonzero slope on $\pa (0,L)$. Indeed, $u_1$, and thus $v_1$ given by \eqref{superposition}, has nonzero slope at $x=0,\, L$ and since $v_1 \leq u$ with the same support $(0,L)$, the same property holds for $u$. 
It follows that there are at least two solutions of \eqref{eq:ob1d*} when $L=kL_c$ for  $k\in \NN $ sufficiently large.
\end{remark}

\subsection{Proof of Theorem \ref{thm:threshold} (for $f$ concave)}

 
 \noindent{\bf Extinction:}
 First, we prove Theorem  \ref{thm:threshold} in the cases where we observe extinction.
 
 We recall that $p_\infty \in L^2(0,T,H^1(\R))$, and so $x\to p_\infty (x,t) $ is in $H^1(\R)$ and therefore continuous for a.e. $t>0$. 
For such $t$, we assume that the open set $\{p_\infty(\cdot,t)>0\}$ is not empty. It can then be written as the union of its connected components $(a_i,b_i)$ and 
Proposition \ref{prop:1} together with the assumption that $f$ is concave implies that (for a.e. $t>0$)
$$ - p_\infty '' = g \leq f(p_\infty) \mbox{ in } (a_i,b_i), \qquad p_\infty>0  \mbox{ in } (a_i,b_i).$$

When $\int_0^{p_M} f(s)\, ds <0$, Proposition \ref{prop:L0} gives a contradiction (see Remark \ref{rmk:L_0}). It follows that
$p_\infty(x,t)=0$ for all $x\in\R$ and a.e. $t>0$.
Equation \eqref{eq:asg} and the concavity of $f$ then implies
$$\pa_t \rho_\infty = \rho_\infty g \leq \rho_\infty f(p_\infty) = \rho_\infty f(0)$$
 and so $\rho_\infty(x,t)\leq \rho_{in} e^{f(0)t}\to 0$ as $t\to \infty$.

When $\int_0^{p_M} f(s)\, ds >0$ and $\rho_{in}=\chi_{(0,L)}$ with $L<L_0$, we will also get a contradiction from Proposition \ref{prop:L0}  if we can show that $|b_i-a_i|<L_0$ (using Remark \ref{rmk:L_0}).
To show that  $|b_i-a_i|<L_0$, we note that 
$$ \frac d {dt} \int_{\R} \rho_m(x,t)\, dx = \int_{\R} \rho_m f(p_m)\, dx   \leq C \int_{\R} \rho_m(x,t)\, dx$$
and we deduce that 
$$ 
\int_{\R} \rho (x,t) \, dx \leq \int_{\R} \rho_{in} (x) \, dx + Ct  < L_0 + Ct \qquad \mbox{ a.e } t>0.
$$
In particular there exists $\eta>0$ such that 
\begin{equation}\label{eq:rhoL}
\int_{\R} \rho (t,x) \, dx < L_0 \qquad \mbox{ a.e. $t\in(0,\eta)$.}
\end{equation}
 In view of Lemma \ref{lem:HS}, we can always assume that
$\int p_\infty(x,t) (1-\rho_\infty(x,t)) dx = 0$ so that $\rho_\infty(x,t) =1$ a.e. in $(a_i,b_i)$ and \eqref{eq:rhoL} gives
$ |b_i-a_i|< L_0 $
(note that this argument only requires $\int_{\R} \rho_{in} (x) \, dx <L_0$).
and Proposition \ref{prop:L0} again gives that $p_\infty(x,t) = 0 $  $\forall x\in \R$, a.e. $t\in(0,\eta)$. As above, it follows that 
$$\pa_t \rho_\infty = \rho_\infty g \leq \rho_\infty f(p_\infty) = \rho_\infty f(0) \mbox{ in } \R\times (0,\eta).$$
and so $\int_{\R} \rho(x,\eta) \, dx \leq \int_{\R} \rho_{in}(x) \, dx$. We can then use the same argument to prove that $p_\infty(\cdot,t) = 0 $  a.e. $t\in(0,2\eta)$.
Successive iterations yield the result.

\bigskip

 \noindent{\bf Invasion:} 
We now prove the last part  of Theorem \ref{thm:threshold}.
Given $L>L_0$, we consider $q(x)$, a solution of  \eqref{eq:ob1d*}, which we extend to $\R$ by $0$. 
Assumption \eqref{eq:pinch} implies that  $p_{in,m}(x)\geq q(x)$, so we can compare $p_m(x,t)$ with the solution of the porous media equation with initial condition $q(x)$, which is better behaved than $p_m(x,t)$ thanks to Proposition \ref{prop:monotone}.
We note that we can always assume that $q'(0) >0$. Indeed, \eqref{eq:ob1d*} has solution with $q'(0)=0$ only for discrete values of $L$ (see Remark  \ref{rem:u0}). If $L$ happens to be one of those values, we can always take a solution $q(x)$ of \eqref{eq:ob1d*} with a smaller value of $L$ (still greater than $L_0$) so that $q'(0)>0$.

\medskip

We recall that $q$ satisfies $\frac 1 2 |q'(x)|^2 + F(q(x)) = F(M)$ with
$$ 
M:=\sup_{x\in(0,L)} q(x)  
$$
and $|q'(0)|=|q'(L)| = \sqrt{2F(M)}>0$.

\medskip

With this function $q(x)$, we have the following lemma, which immediately implies 
the last part  of Theorem \ref{thm:threshold}:
 
\begin{lemma}
Let $\tilde{\rho}_m$ be the solution of \eqref{rho_m} with $\tilde p_m (x,0) = q(x)$. Then 
\item[(i)] For all $m$ the functions $\tilde{\rho}_m$ and $\tilde{p}_m$ are monotone increasing in time.
\item[(ii)] The limit $(\tilde \rho_\infty,\tilde p_\infty)$ of $(\tilde{\rho}_m,\tilde{p}_m)$ along any convergent subsequence satisfies
\begin{equation}\label{eq:limit1D}
\pa_t\tilde  \rho_\infty = \pa_{xx} \tilde p_\infty + \tilde \rho_\infty f(\tilde p_\infty) , \qquad \tilde p_\infty \in P_\infty(\tilde \rho_\infty). 
\end{equation}
\item[(iii)] For all $t>0$, we have $\tilde \rho_\infty(\cdot,t) = \chi_{\Omega(t)}$ for some open set $\Omega(t)$ satisfying 
$$ \Omega(t) \supset (- \sqrt{2F(M)}t,L+\sqrt{2F(M)}t).$$
\end{lemma}

\begin{proof}

 (i) and (ii)  immediately follow from Proposition \ref{prop:monotone} since $q$ satisfies
$$ q''+f(q) \chi_{\{q>0\}} \geq 0 \qquad \mbox{ in } \R.$$
Furthermore, since the functions are monotone increasing in time, we can use the argument in \cite{MPQ}  to prove that $\tilde \rho_\infty$ is a characteristic function:
let $w(x,t) := \int_0^t e^{-f(0)s} \tilde p_\infty(x,s)\, ds$, which solves (using \eqref{eq:limit1D})
$$ \pa_{xx} w(t) = e^{-f(0)t}\tilde \rho_\infty(t) -\tilde  \rho_{in} + \int_0^t e^{-f(0)s} \tilde \rho_\infty(s) [f(0)-f(\tilde p_\infty(s)]\, ds \mbox{ a.e. $x$  in $\R$. }$$
Since the right hand side is bounded in $L^\infty(\R)$, we see that for all $t>0$, $w(\cdot,t) \in W^{2,\infty}(\R)$. As a consequence,  
the set  $\Omega(t) := \{w(t)>0\} $ is open and  $\pa_{xx} w=0$ a.e. in $\{w(t)=0\}=\R\setminus\Omega(t)$  (see Remark \ref{rem:11} below).
Finally we check that $\tilde \rho_{\infty}(\cdot,t)=\chi_{\Omega(t)}$. Since $\Omega(t)\subset \{\tilde p_\infty(t)>0\}$, we have $\tilde \rho_\infty(x,t)=1$ a.e. in $\Omega(t)$.
On the other hand, for $x\in \R\setminus\Omega(t)$, we have $\tilde p_\infty(x,s)=0$ a.e. $s\in(0,t)$ by definition of $w$ and thus
$$0 = e^{-f(0)t}\tilde \rho_\infty(t) -\tilde  \rho_{in} \mbox{ a.e. in } \R\setminus\Omega(t).$$
Due to the choice of initial data we have $\tilde \rho_{in}=0$ in $\{w(t)=0\} \subset \{\tilde p_{in}=0\}$, and so it follows that $\tilde \rho_\infty(x,t) = 0 $ a.e. in $ \R\setminus\Omega(t)$. 
Summarizing, we showed that 
\begin{equation}\label{eq:rho1d} \tilde \rho_\infty(x,t) = \chi_{\Omega(t)}, \quad \Omega(t) := \{w(t)>0\} \end{equation}
From the Hele-Shaw condition we have $\tilde p_{\infty}(t)=0$ a.e.  in $ \R\setminus\Omega(t)$, and thus we also have 
$$|\Omega(t) \Delta \{ \tilde p_\infty(t)>0\}| =0.$$

\medskip

Equation \eqref{eq:limit1D}, together with  \eqref{eq:rho1d} gives, for all smooth test function $\phi(x)$:
\begin{align*}
 \frac{d}{dt} \int_{\Omega(t)} \phi(x)\, dx = \frac{d}{dt} \int_\R \tilde \rho_\infty (x,t)\phi(x)\, dx
 & =  \int_\R \pa_t \tilde \rho_\infty(x,t) \phi(x)\, dx\\
& = \int_\R \tilde p_\infty \pa_{xx}\phi + \tilde \rho_\infty f(\tilde p_\infty) \phi\, dx \\
& = \int_{\Omega(t)}\tilde  p_\infty \pa_{xx}\phi +  f(\tilde p_\infty) \phi \, dx .
\end{align*}

As mentioned before, since we are in one space dimension. Proposition \ref{prop:monotone} implies that $ \tilde p_\infty(\cdot,t)$ solves $-\pa_{xx} \tilde p_\infty = f(\tilde p_\infty)$ on any connected components of $\Omega(t)$ and so
\begin{equation}\label{eq:omega}
\int_{\Omega(t)} \tilde p_\infty \pa_{xx}\phi +  f(\tilde p_\infty) \phi \, dx  = \int_{\pa\Omega(t)} |\pa_x \tilde p_{\infty}| \phi
\end{equation}
Furthermore, if $(a(t),b(t))$ is the connected component of $\{\tilde p_{\infty}(\cdot,t)>0\}$ containing $(0,L)$ (recall that $\tilde p_{\infty}(x,t)\geq q(x)$ and so $\tilde p_{\infty}(\cdot ,t)>0$ in $(0,L)$ for all $t\geq 0$), we can write 
$$\frac 1 2 |\pa_x \tilde p_\infty(x)|^2 + F(\tilde p_\infty(x)) = F\left(\sup_{(a,b)} \tilde p_\infty\right)\geq  F\left(\sup_{(a,b)} q \right) =F(M) \mbox{ for all } x\in[a,b]$$
and so
$$|\pa_x \tilde p_\infty(a(t)) | = |\pa_x \tilde p_\infty(b(t)) | \geq \sqrt{2F(M)}.$$
We thus have
$$
\frac{d}{dt} \int_{\Omega(t)} \phi(x)\, dx 
\geq \sqrt{2F(M)} [\phi(a(t))+\phi(b(t)] \qquad \mbox{ in } \mathcal D'(0,\infty)$$
which implies 
$$
\int_{\pa \Omega(t)} V \phi(x) d\mathcal H^{0}(x) \geq \int_{\pa (a(t),b(t))} \sqrt{2F(M)}  \phi(x) d\mathcal H^{0}(x)
$$
where  $V$ denotes the normal velocity of $\pa \Omega(t)$ (which is a countable set of points).
Since this holds for any test function $\phi$, we deduce that $V\geq 0$ on $\pa\Omega(t)$ and $V\geq  \sqrt{2F(M)} $ on $\pa (a(t),b(t))$.
The result follows.


\end{proof}

\begin{remark}\label{rem:11}
We recall that if $u\in W^{1,1}_{loc}(\Omega)$ and  if we denote $E_\alpha= \{u=\alpha\}$ for any $\alpha\in \R$, then $\na u(x) =0$ for almost every $x\in E_\alpha$.
In our setting, this classical result implies that $\pa_x w = 0 $ a.e. in $\{w=0\}$ and that $\pa_{xx} w =0$ a.e. in $\{\pa_x w=0\}$. The fact that $\pa_{xx} w =0$ a.e. in $\{w=0\}$ follows.
Importantly, this argument does not require any particular properties for the boundary $\pa \{w>0\}$.
\end{remark}

\subsection{Proof of Proposition \ref{prop:long} (general $f$)}
We now assume that 
$ p_{in,m}(x) \leq p_M \chi_{(0,L)}(x)$
with $L$ such that
\begin{equation} \label{eq:L*}
L<L^* = \frac{\pi}{\sqrt K }.
\end{equation}
(with $K$ defined by \eqref{eq:K}).

The proof is similar to the first part of the proof of Theorem  \ref{thm:threshold}.
First we remark that \eqref{eq:K} implies $f(p_m) \leq K p_m,$
and so  we have $g \leq K p_{\infty}$, with the notations of Proposition \ref{prop:1}.

Recall that $p_{\infty}(\cdot,t)$ is continuous for a.e. $t>0$. For such $t$, we assume that the open set $\{p_\infty(\cdot,t)>0\}$ is not empty. It can then be written as the union of its connected components $(a_i,b_i)$ and 
Proposition \ref{prop:1} together with the remark above implies that, for a.e. $t>0$, $p_{\infty}(\cdot,t)$ solves
$$ - p_\infty '' = g \leq  K p_\infty \mbox{ in } (a_i,b_i), \qquad p_\infty>0  \mbox{ in } (a_i,b_i).$$

We can now show (as in the proof Theorem  \ref{thm:threshold} in the extinction case) that
$ |b_i-a_i | < L^*$ for some small time $t\in(0,\eta)$. Proceeding as in the last part of proof for  Proposition \ref{prop:L0}, we see that Poincar\'{e} inequality yields $p_\infty(x,t) = 0$ in $(a_i,b_i)$.
Since this holds for any connected components of $\{p_\infty(\cdot,t)>0\}$ it follows that $p_\infty(x,t) = 0$  for a.e. $t\in(0,\eta)$.

This in turn implies that $\rho_{\infty}(x,t) = \rho_{in}(x) e^{f(0)t}$ for $t<\eta$, so that the support of $p_\infty(x,\cdot)$ has not grown and we can iterate to show that $p_\infty(x,t)=0$ for all $t>0$.

\section{Traveling wave solutions}
First we study the existence of traveling waves for the limiting problem  and prove Proposition \ref{prop:limtw}).
We start with the pressure equation:
Recalling that 
$F(p)= \int_0^p f(s)\,ds$, we have:

\begin{proposition}\label{prop:h}

\item[(i)] if $f$ satisfies  \eqref{eq:f1'} (monostable), then
the equation
\begin{equation}\label{eq:heq}
h'' + f(h) = 0, \qquad h\geq 0 \mbox{ in } (-\infty,0), \quad h(0)=0
\end{equation}
has a unique nonzero, bounded solution which is decreasing. This solution satisfies $h'(0)= -\sqrt{2F(p_M)}$ and $\lim_{x\to -\infty} h = p_M$. 
\item[(ii)] If  $f$ satisfies  \eqref{eq:f1} (bistable) then
\begin{itemize} 
\item  if $F(p_M)>0$, then
the equation \eqref{eq:heq} has two bounded solutions:
One monotone solution as in (i) and one periodic solution with $h'(0)=0$.
\item if $F(p_M)=0$, then \eqref{eq:heq} has a unique bounded solution which is  decreasing. This solution  satisfies   $h'(0)=0$ and $\lim_{-\infty} h = p_M$.
\item
 if $F(p_M)<0$, then \eqref{eq:heq} has no bounded solution.
\end{itemize}
\end{proposition}


\begin{proof}[Proof of Proposition \ref{prop:h}]
Assume that $h$ is a bounded solution of  \eqref{eq:heq}.

First, we note that  if  there is $x_0<0$  such that $h'(x_0)=0$ then uniqueness principle for second order ODE on $[0,x_0]$ with $h(x_0)$ and $h'(x_0)=0$ fixed yields that $h(x_0+s) = h(x_0-s)$ for all $s\in[0,x_0]$. This implies that $h(2x_0)=h(0)=0$ and $h'(2x_0) = - h'(0)$.
Since $h\geq 0$ in $(-\infty,0)$, we must then have $h'(0)=0$. From the uniqueness principle again, it follows that $h(2x_0+s) = h(s)$, namely $h$ is periodic.
When $h'(0)\neq 0$, we deduce that $h'$ cannot vanish in $(-\infty,0)$ and is thus always negative.

Next we show that $h <p_M$. Note that $h''\geq -f(h)> 0$ if $h > p_M$, and thus $h$ is strictly convex when it is above $p_M$. Hence if $h(x_1)>p_M$ for some $x_1<0$, then $h$ must be unbounded in $(-\infty,0)$, a contradiction. If $h(x_1)=p_M$, then either $h$ is unbounded or $h'(x_1)=0$, but in this last case, we would then have $h(x)=p_M$ for all $x$ which contradicts the fact that $h(0)=0$.

\medskip

We thus have that $h$ is either periodic or monotone decreasing and that $0\leq h<p_M$ in $(-\infty,0)$. 
When $h$ is nontrivial and monotone, $\lim_{x\to -\infty} h(x) =\ell \in (0, p_M]$  exists and satisfies $f(\ell)=0$, that is $\ell = \alpha$ or $p_M$.  

\medskip
Next let us show that $h$ can be nontrivial and periodic only when $F(p_M)>0$ and $f(0)<0$. Multiplying \eqref{eq:heq} by $h'$ and integrating, we find
\begin{equation}\label{eq:constant}
\frac{1}{2}|h'|^2 +F(h) =   \frac{1}{2}|h'(0)|^2  \hbox{ in } (-\infty,0).
\end{equation}
When $h$ is periodic, we have $h'(0)=h'(x_0)=0$. Hence it follows that $F(h(z_0))=0$ for some $h(z_0)\in (0,p_M)$. This situation can only happen when $f(0)<0$ and $F(p_M)>0$.

\medskip

Lastly, when $h$ is monotone, \eqref{eq:constant} implies $F(\ell) =  \frac{1}{2}|h'(0)|^2 $.
If $F(p_M)<0$, then this is impossible and no bounded solution exist, thus proving the last part of the proposition.
If $F(p_M)\geq0$, and since $F(\alpha)<0$, we deduce $\ell = p_M$ or $\ell=0$. The latter case yields a trivial $h$ when $f(0)=0$. In the former case we have 
$$ h'(0) = - \sqrt{2 F(p_M)}.$$
The uniqueness of $h$ follows since the initial value problem $h'' + f(h) = 0$ with $h(0)=0$ and $h'(0) = -\sqrt{2 F(p_M)}$ has a unique solution. 
\medskip

The existence of $h$  when $F(p_M)\geq 0$ can now be obtained by solving
$$
h''+f(h)=0 \hbox{ in } (-\infty,0), \quad h(0)=0 \hbox{ and } h'(0)= -\sqrt{2F(p_M)} . $$

To check that $h$ is bounded, we will use \eqref{eq:constant}. Observe that $h <p_M$: if $h(x_0) = p_M$  for some $x_0<0$ then $h'(x_0)=0$, and the uniqueness principle for the ODE yields $h(x) \equiv p_M$ for all $x<0$, a contradiction. From \eqref{eq:constant} 
and the fact that $F$ has a maximum at $p_M$ in $[0,p_M]$ it follows that $h'(x)$ cannot vanish, and so $h$ is monotone decreasing and bounded above and therefore positive and bounded in $(-\infty,0)$.
\end{proof}

\medskip

We can now prove the existence and uniqueness of traveling wave for the limiting equation:
\begin{proof}[Proof of Proposition \ref{prop:limtw}]
We start with the following remark:
Since $\bar \rho$ is non-increasing and $\bar \rho(+\infty)=\ell<1$, there must exist $a \in [-\infty,+\infty)$ such that $\bar \rho(x)=1$ for $x\in (-\infty,a)$ and $\bar \rho(x)<1$ in $(a,+\infty)$. 
If $a=-\infty$, then $\bar \rho<1$ and $\bar p=0$ in $\R$ 
which contradicts \eqref{eq:TWbc}.
We can thus always assume that $a\in(-\infty,+\infty)$ and up to translation we will take $a=0$.

Let us first classify possible pressure profiles. Up to a translation, any traveling wave pressure must satisfy 
 $\bar p=0$ in $(0,+\infty)$ and
$ \bar p''+f(\bar p) = 0$ in $(-\infty,0)$. 
In particular, $\bar p$ coincides on $(-\infty,0)$ with a  monotone decreasing  solution of \eqref{eq:heq}.
Using Proposition \ref{prop:h}, we deduce that there is no such $\bar p$ when $F(p_M)<0$ which implies (i).  On the other hand, if  $F(p_M)\geq 0$ then again Proposition \ref{prop:h} yields that $\bar p$ is uniquely determined and satisfies
 $\bar p'(0)= -\sqrt{2F(p_M)}$. 
 

\medskip

Next we consider possible travelling wave densities. To this end, we observe that \eqref{eq:limtw} is equivalent to
\begin{equation}\label{eq:delta}
 -c\bar \rho'  = - \bar p'(0) \delta_{x=0} + \bar \rho f(0) \chi_{(0,\infty)}= \sqrt{2F(p_M)} \delta_{x=0}+ \bar \rho f(0) \chi_{(0,\infty)}.
\end{equation}


\medskip

When $F(p_M)=0$ and $f(0)<0$,  \eqref{eq:delta} yields a continuous solution
$$
\bar \rho(x) = \chi_{(-\infty,0)} (x)+ e^{-\frac{f(0)}{c} x} \chi_{(0,+\infty)}(x),
$$
which satisfies the boundary conditions with $\ell=0$ for all $ c<0$ and 
$$
\bar \rho(x) = \chi_{(-\infty,0)} (x)
$$
with $c=0$. This implies (iii).

\medskip

When $F(p_M) >0$, $c$ is positive due to \eqref{eq:delta} at $x=0$ with our assumption  $\bar \rho'\leq 0$. So \eqref{eq:delta} again gives $\bar \rho f(0) \geq 0$ in $(0,\infty)$. 
When $f(0)<0$, it follows that $\bar \rho=0$ in $(0,+\infty)$, yielding $\ell=0$, concluding (ii). 

\medskip

Lastly  let us consider the monostable case $f(0)=0$. In this case \eqref{eq:delta} implies  $-c\bar \rho' = 0 $ in $(0,+\infty)$, which yields $\bar \rho=\ell$  in $(0,+\infty)$. Since $\bar{\rho}=1$ in $(-\infty,0)$, it follows from \eqref{eq:delta} that $c (1-\ell) = \sqrt{2F(p_M)}$. We now conclude.


\medskip

 The existence of the traveling wave when $F(p_M)> 0$ is proved by taking $\bar p=h$ as the unique monotone solution of \eqref{eq:heq}, $c= \frac{\sqrt{2F(p_M)}} {1-\ell}$ and  $\rho = \chi_{(-\infty,0)}+\ell \chi_{(0,+\infty)}$ and checking that this solves \eqref{eq:limtw} for appropriate values of $\ell$.

\end{proof}

\begin{proof}[Proof of  Theorem \ref{thm:TW} (i) (Receding Traveling Waves)]
Let us fix $m>1$ and denote by $(\bar \rho (x), c) $ the traveling wave of \eqref{rho_m} given by Theorem \ref{thm:TW0}, and by $\bar p(x)$ the corresponding pressure given by $\bar p(x) = \frac m {m-1} \bar \rho(x)^{m-1}$. Since we assume that $F(p) = \int_0^p f(u)\, du$ satisfies $F(p_M)<0$, we have $c<0$ for large $m$.

\medskip

Recall from Theorem~\ref{thm:TW0} that $\bar{\rho}$ and thus $\bar{p}$ is monotone decreasing and that $0\leq \bar p(x)\leq p_M$. Furthermore, we claim that there exists a constant $C>0$ independent of $m$ such that 
\begin{equation}\label{eq:p'p}
\bar p'(x)^2 \leq C \bar p(x) \leq Cp_M.
\end{equation}
Indeed, since $\pa_t (\bar \rho(x-ct))= -c \bar{\rho}'\leq 0$, \eqref{rho_m}  gives
$ \bar \rho'\bar p'+ \bar \rho \bar p'' \leq -\bar  \rho f(p) $. Since $\rho' p'\geq 0$ and $\bar\rho$ is positive, $p''(x)\leq \sup_{0\leq p \leq p_M} -f(p) \leq C$. Multiplying by $\bar p'(x)$ and integrating over $(x,\infty)$, we find  \eqref{eq:p'p}.

\medskip

Next, set $\eta := -F(p_M)/(2p_M)$ and consider the function
$$ g(x)= \frac 1 2 \bar p'(x)^2 + F(\bar p(x)) + \eta \bar  p(x).$$
We claim that there exists $\bar x\in \R$ such that 
\begin{equation}\label{eq:xbar}
\begin{cases}
g(\bar x) \leq  \frac{F(p_M)}{4} \\
g'(\bar x) \geq 0 \\
\bar p(\bar x) \geq \eta_1 
\end{cases}
\end{equation}
for some small $\eta_1>0$.
Indeed, we have $g(-\infty) = F(p_M) + \eta p_M = \frac{F(p_M)}{2} <0$ and $g(+\infty) =0$, which implies the existence of $\bar x$ satisfying the first two conditions. The first condition implies in particular
$$   F(\bar p(\bar x)) \leq  \frac{F(p_M)}{4} -\frac 1 2 \bar p'(\bar x)^2 - \eta \bar p(\bar x)  \leq  \frac{F(p_M)}{4}<0$$
which gives the third condition in \eqref{eq:xbar}.
Since $g'(x) = \bar p'(x) ( \bar p''(x) + f(\bar p(x)) + \eta) $, the condition $g'(\bar x) \geq 0$ implies 
$$ \bar p''(\bar x) + f(\bar p(\bar x)) \leq -\eta$$
and so
$$ \bar p(\bar x) (\bar p''(\bar x) + f(\bar p(\bar x))) \leq -\eta \eta_1.$$
The pressure equation now gives
$$ -c\bar p'(\bar x) = (m-1) \bar p(\bar x)  (\bar p''(\bar x)  +f(\bar p(\bar x) )) +\bar  p'(\bar x) ^2 \leq -(m-1) \eta \eta_1 +\bar  p'(\bar x) ^2$$
that is (recall that $0\leq - \bar p'\leq C$)
$$ c  \leq -(m-1)\frac{ \eta \eta_1}{C} +C$$
and \eqref{eq:cm0} follows.

\end{proof}

\begin{proof}[Proof of  Theorem \ref{thm:TW} (ii) (Advancing traveling Waves)]
We now assume that $F(p_M)=\int_0^{p_M} f(p)\, dp >0$. As before we denote by $\bar \rho_m$, $\bar p_m$, $c^*_m>0$ the traveling wave given by Theorems \ref{thm:TW0} and \ref{thm:TW0'} and its associated pressure and velocity for fixed $m$. For $m$ large enough, we have $ \int_0^{\rho_m^+} \rho^m f\left(p_m)\right)\,d\rho >0$.
We also denote by $\bar \rho^*, \bar p^* ,c^*>0$ the unique traveling wave of the limiting problem given by Proposition \ref{prop:limtw}. 
Finally, we write $p_m(x,t) :=\bar p_m(x-c^*_mt)$.

\medskip

Since we have uniqueness of the traveling waves for the limiting problem, it will be easy to conclude if we can show that the limit of $\bar \rho_m,\bar p_m$ solves \eqref{eq:limtw} and are supported in $(-\infty,0)$.
In order to  get   \eqref{eq:limtw}, the main step is to show that $c^*_m$ converges which requires an upper bound on $c^*_m$. To show that $\bar \rho_m,\bar p_m$ do not degenerate and that their limits are positive in $(-\infty,0)$, we will need to know that $c_m^*$ is uniformly bounded away from $0$ when $m$ is large.
\medskip

We thus start the proof by proving that $c^*_m\to c^*$, by constructing appropriate barriers for the pressure equation.
We recall the pressure equation, solved by $p_m(x,t)$:
\begin{equation}\label{eq:pressure eq}
p_t = (m-1) p(\pa_{xx}p + f(p)) + |\pa_x p|^2.
\end{equation}
Set $f_\eps(p) := f(p)+\e$. Since $f(p)<0$ for $p>p_M$, if $\eps$ is small then there exists $p_M^\eps > p_M $ with $p_M^\eps = p_M+ o(1)$ such that $f_\eps(p)<0$ for $p>p_M^\eps$.
Let $h_\eps $ be the monotone decreasing solution of \eqref{lim_tw} with $f_\eps $ instead of $f$, that is:
$$
\begin{cases}
h_\eps''(x) = f_\eps(h_\eps(x) ) \quad \mbox{ for } x\in (-\infty,0), \\
h_\eps(-\infty) = p_M^\eps, \quad h_\eps(x) = 0 \mbox{ in } [0,\infty)
\end{cases}
$$
(we can proceed as in Proposition \ref{prop:h} to show that such an $h_\eps$ exists). 
We then claim that the function $p_\eps(x,t) = h_\eps (x-c_\e t)$ with $c_\e = |h_\eps'|(0)+ \e$ is a supersolution
for \eqref{eq:pressure eq}.
Indeed, we have
$$
p_\eps(\pa_{xx}p_\eps + f(p_\eps))   = p_\eps( -f_\eps (p_\eps) + f(p_\eps)) =  -\e p_\eps  \qquad \mbox{ in }\{p_\eps>0\} = \{x-c_\e t <0\}
$$
and
\begin{align*}
\pa_t p_\eps (x,t)- |\pa_x p_\eps|^2 & = c_\eps |h_\eps'|(x-c_\e t) - |h_\eps'|^2 (x-c_\e t) \\
& = (|h_\eps'|(0)| + \e) |h_\eps'|(x-c_\e t) - |h_\eps'|^2 (x-c_\e t) \\
& \geq  0  \qquad \mbox{ when  } -\delta \leq  x-c_\e t \leq 0 
\end{align*}
for some small $\delta$.
We thus have
$$ 
\pa_t p_\eps - (m-1) p_\eps(\pa_{xx} p_\eps + f(p_\eps)) - |\pa_x p_\eps|^2 
\geq  
\begin{cases}
0 & \mbox{ if }   x-c_\e t \geq -\delta\\
c_\eps\pa_x p_\eps - |\pa_x p_\eps|^2 + \eps (m-1) p_\eps
& \mbox{ if }   x-c_\e t \leq -\delta
\end{cases}
$$
Since $p_\eps$ is bounded below when  $ x-c_\eps t < - \delta$, this last term is non-negative of $m$ is large enough. 

\medskip

Having verified the supersolution property of $p_{\eps}$, let us now compare $p_m$ with $p_\eps$. Let now $x_0$ be such that $p_\eps(x,0) \geq p_M$ for $x \leq -x_0$.
Since $p_m(x,0)=0$ for $x>0$ (see Theorem \ref{thm:TW0}), we have $p_m(x,0) \leq p_\eps (x + x_0,0)$ and the comparison principle implies that $p_m(x,t) \leq p_\eps (x + x_0,t)$  for all $t>0$.
We deduce that
$$ c_m  \leq c_\eps $$ 
for sufficiently large $m$ and so
$ \limsup_{m\to \infty} c_m \leq c_\eps$.
Furthermore, it is easy to check that
$$
c_\e = \sqrt{2F_\eps(p^\eps_M)} + \e \to c^* \mbox{ as } \eps\to 0
$$
hence
$$ \limsup_{m\to \infty} c_m \leq c^*. $$



Similarly, we can construct a subsolution of \eqref{eq:pressure eq}  by using the function  $f_\eps(p) = f(p)-\e$ and $c_\eps =  |h_\eps'|(0)- \e$.
We then have:
 $$
\pa_t p_\eps - (m-1) p_\eps(\pa_{xx} p_\eps + f(p_\eps)) - |\pa_x p_\eps|^2 
=(|h_\eps'|(0) - \e) |h_\eps'|(x-c_\e t) - |h_\eps'|^2 (x-c_\e t) - \eps (m-1) h_\eps (x-c_\e t) 
$$
which is negative when $-\delta\leq x-c_\e t \leq 0$ if $\delta$ is such that $|h_\eps'(0)| - \e \leq |h_\eps'(y)| $ for $y\in (-\delta,0)$ and negative for $  x-c_\e t \leq -\delta$, provided $m$ is large enough depending on $\eps$ and $\delta$.
Proceeding as above, we deduce that:
$$
\liminf c_m \geq c^*.
$$
\medskip

We have thus showed that $ c_m \to c^*$ and it remains to show that $\bar \rho_m \to \chi_{(0,\infty)}$ and $\bar p_m\to h$, the monotone solution of \eqref{eq:heq}.

\medskip

Since $\bar \rho_m$ and $\bar p_m$ are monotone decreasing and uniformly bounded, they are uniformly bounded in $BV(\R)$. Hence along a subsequence $\bar\rho_m$ and $\bar p$ converge  to monotone functions $\bar \rho_\infty$ and $\bar p_\infty$ strongly in $L^1_{loc}$ and almost everywhere.
Proceeding as in  the proof of Proposition  \ref{prop:monotone} we get that $\bar \rho_\infty$ and $\bar p_\infty$ solve \eqref{eq:limtw} (with $c=c^*$).
In order to conclude, it remains to show that $\bar \rho_\infty$ and $\bar p_\infty$ satisfy\eqref{eq:TWbc}.

\medskip

Let us first show that $p_\infty(x)>0$ in $(-\infty,0)$. We use the fact that $c^*_m \to c^*>0$ and so $c^*_m \geq  c^*/2>0 $ for $m$ large enough and that $\bar  p_m(x)$ solves
\begin{equation}\label{first}
c^*_{m} (- \bar p_m') - | \bar p_m'|^2 = (m-1)\bar p_m(\bar p_m'' + f(\bar p_m)) 
\end{equation}
In particular, we must have $(\bar p_m)''+f(\bar p_m)\geq 0$ whenever $\bar p_m' \in (-c^*_m,0)$, and thus $\bar p_m'' >-f(\bar p_m) \geq -K\bar p_m$ for this range of $\bar p_m'$. 
Finally, we recall that  $\bar p_m(0)=0$ with $\bar p'_m(0^-) =-c^*_m < c^*/2$.
It follows that $\bar p_m' < -c^*/4$ on $[-\eta,0]$, for some $\eta$ independent of $m$ and thus  $\bar p_m(x) > \frac{c^* x}{4}$ for $x\in [-\eta,0)$. From the monotonicity of $\bar p_\infty$ implies that $\bar p_\infty$ is positive in $(-\infty,0)$.

\medskip

Since $\mathrm{Supp } \, \rho_m=(-\infty,0)$, we have $\bar \rho_\infty \leq \chi_{(0,\infty)}$ and $\bar p_{\infty}=0$ on $(0, \infty)$. Since $\bar p_\infty(1-\bar \rho_\infty)=0$, we deduce that $\bar \rho_\infty =\chi_{(0,\infty)}$. Now Proposition \ref{prop:h} implies $\bar p_\infty=h$, and we can conclude.

\end{proof}

\bibliography{DrugTumor_ref}
\bibliographystyle{plain}
\end{document}